\documentclass[11pt]{article}
\usepackage{amssymb,latexsym}
\usepackage{amsmath}

\title{ \large {\bf Identities of
finitely generated graded algebras with involution.}}

\author{ Irina Sviridova\thanks{e-mail \texttt{I.Sviridova@mat.unb.br}}
\\
\\
Departamento de Matem\'atica,\\
Universidade de Bras\'\i lia,\\
70910-900 Bras\'\i lia, DF, Brazil }

\date{October 09, 2014.}

\newtheorem{theorem}{Theorem}[section]
\newtheorem{lemma}{Lemma}[section]

\newtheorem{corollary}[lemma]{Corollary}
\newtheorem{definition}[lemma]{Definition}

\newtheorem{conjecture}{Conjecture}[section]
\newtheorem{assumption}{Assumption}[section]

\begin{document}
\maketitle

\begin{abstract}
We consider associative algebras with involution graded
by a finite abelian group $G$ over a field of characteristic zero.
Suppose that the involution is compatible with the grading.
We represent conditions permitting PI-representability of such algebras.
Particularly, it is proved that a
finitely generated $(\mathbb{Z}/q \mathbb{Z})$-graded
associative PI-algebra with involution
satisfies exactly the same
graded identities with involution as some
finite dimensional $(\mathbb{Z}/q \mathbb{Z})$-graded algebra
with involution for any prime $q$ or $q=4.$
This is an analogue of the theorem of  A.Kemer for ordinary
identities \cite{Kem1}, and an extension of the result
of the author for identities with involution \cite{Svi2}.
The similar results were proved also
for graded identities \cite{AB}, \cite{Svi1}.

\bigskip


\end{abstract}

\section*{Introduction}

One of the central problem of the theory of varieties and
PI-theory is the Specht problem: the problem of existence of a finite base for
any system of identities. The original Specht problem \cite{Specht} was formulated for identities of associative algebras over a field of characteristic zero. It was solved positively by Alexander Kemer
\cite{Kem1}, \cite{Kem3}, \cite{Kem5}. The principal and the most difficult part of
the Kemer's solution was the proof of PI-representability of finitely generated PI-superalgebras
\cite{Kem1}, \cite{Kem4}. An algebra (a superalgebra) is called {\it PI-representable} if it satisfies
the same identities ($(\mathbb{Z}/2 \mathbb{Z})$-graded identities) as some finite dimensional
algebra (superalgebra).
The phenomena of PI-representability of finitely generated algebras has also the proper interest.
It is an intriguing question, what are the classes of algebras and identities such that their finitely
generated algebras satisfy the same identities as finite dimensional algebras.

The first results on PI-representability of associative algebras belong to A.Kemer.
He proves that any finitely generated associative $(\mathbb{Z}/2 \mathbb{Z})$-graded PI-algebra
over a field of characteristic zero satisfies the same $(\mathbb{Z}/2 \mathbb{Z})$-graded identities
as some finite dimensional $(\mathbb{Z}/2 \mathbb{Z})$-graded algebra over the same field \cite{Kem1}, \cite{Kem4}. Later, he proves also that a finitely generated associative PI-algebra
over an infinite field satisfies the same ordinary (non-graded) polynomial identities
as some finite dimensional algebra \cite{Kem2}. PI-representability of finitely generated
associative algebras over a commutative associative Noetherian ring with respect to ordinary polynomial
identities was studied in \cite{B1}-\cite{B9}.

The series of results were obtained also for graded identities and
identities with involution of associative algebras over a field of characteristic zero.
If $G$ is a finite abelian group then a finitely generated $G$-graded associative PI-algebra
over an algebraically closed field of characteristic zero satisfies the same graded
identities as some finite dimensional $G$-graded algebra over the same field \cite{Svi1}.
For more general case of a finite (not necessarily abelian) group $G$
it was proved that a finitely generated $G$-graded associative PI-algebra
over a field of characteristic zero satisfies the same graded
identities as some finite dimensional $G$-graded algebra over some extension
of the base field \cite{AB}. As the direct consequences of \cite{AB}, \cite{Svi1}
we have also the similar results for $G$-identities if $G$ is
a finite abelian group of automorphisms of an associative algebra.

Recently, PI-representability was proved also for identities with involution \cite{Svi2}.
A finitely generated associative PI-algebra with involution
over a field of characteristic zero satisfies the same
identities with involution as some finite dimensional algebra with involution
over the same field.

The special interest to graded identities in the case of characteristic zero is explained
by the super-trick and relation with the Specht problem \cite{Kem1}.
Therefore the problem of PI-representability of graded algebras with involution
is of current interest also.

We consider associative algebras over a field $F$ of
characteristic zero. Further they will be called algebras.

An $F$-algebra $A$ is {\it graded by a group $G$}
({\it $G$-graded algebra}) if $A$ can be decomposed as a direct sum $A
=\bigoplus_{\theta \in G} A_{\theta}$ of its vector subspaces $A_{\theta}$ ($\theta \in G$),
where $A_{\theta} A_{\xi} \subseteq A_{\theta \xi}$ holds for all
$\theta, \xi \in G.$ A subspace $V \subseteq A$ is called {\it graded}
if $V=\bigoplus_{\theta \in G} (V \cap A_\theta).$
We consider only gradings by a finite abelian group.

Anti-automorphism $*$ of the second order of an algebra $A$
over $F$ is called {\it involution}. Algebra with involution is also called
$*$-algebra. An element $a$ of a $*$-algebra $A$
is called symmetric if $a^*=a,$ and
skew-symmetric if $a^*=-a.$ Particularly, $a+a^*$ is symmetric and
$a-a^*$ is skew-symmetric for any element $a \in A.$ It is
clear that $A=A^{+} \oplus A^{-},$ where $A^{+}$ is the subspace
formed by all symmetric elements ({\it symmetric part}), and
$A^{-}$ the subspace of all skew-symmetric elements of $A$
({\it skew-symmetric part}).
We also use the notations $a \circ b =a b + b a,$ and $[a,b]=a b - b a.$
It is clear that the symmetric part $A^{+}$ of
a $*$-algebra $A$ with the operation $\circ$ is a Jordan algebra, and
the skew-symmetric part $A^{-}$ with the operation $[,]$ is a Lie algebra.

Let $G=\{ \hat{\theta}_1=\mathfrak{e},
\hat{\theta}_2,\dots,\hat{\theta}_{\mathfrak{m}} \}$
be a finite abelian group of order $\mathfrak{m}$ with the unit $\mathfrak{e}.$
Let us consider a $G$-graded algebra
$A=\bigoplus_{\theta \in G} A_{\theta}$ with involution.
We assume that involution $*$ is graded anti-automorphism of $A,$ i.e.
$A_{\theta}^{*}=A_{\theta}$ for any $\theta \in G.$ This is equivalent to condition
(see, e.g., \cite{BahtShestZ}) that the subspaces
$A^{+},$ $A^{-}$ are graded. Particularly, we have
$A=\bigoplus_{\theta \in G} (A_{\theta}^{+} \oplus A_{\theta}^{-}),$ where
$A^{\delta}=\bigoplus_{\theta \in G} A_{\theta}^{\delta},$
\ ($\delta \in \{ +, -\}$); and
$A_{\theta}=A_{\theta}^{+} \oplus A_{\theta}^{-},$ \ ($\theta \in G$).
We say that an element $a \in A_{\theta}^{\delta}$ ($\delta \in \{ +, -\},$ \
$\theta \in G$) is {\it homogeneous} of {\it complete degree }
$\deg_{\widehat{G}} a=(\delta,\theta)$ or simply {\it $\widehat{G}$-homogeneous}.

Note that if the base field $F$ contains a primitive root of unity
$\sqrt[\mathfrak{m}]{1}$ of order $\mathfrak{m}=|G|$
then a $G$-grading on $A$ naturally corresponds
to the action on $A$ of the group ${\rm Irr}{G}=\{ \chi_{k} | k=1,\dots,|G| \} \cong G$
of irreducible characters of $G.$ An irreducible character $\chi \in {\rm Irr}{G}$
acts on $A$ as an automorphism associating to any element $a=\sum_{\theta \in G} a_{\theta} \in A$
the element $\chi(a)=\sum_{\theta \in G} \chi(\theta) a_{\theta}$ (\cite{GZbook}).
Then the involution $*$ of $A$ is graded iff it commutes with any $\chi \in {\rm Irr}{G}.$
Thus in this case the group $\widehat{G}={\rm Irr}{G} \times \langle *
\rangle \cong  G \times \mathbb{Z}/2\mathbb{Z}$ of automorphisms and anti-automorphisms acts on $A.$
We refer readers for more details about connection of gradings
with automorphism group actions to \cite{BahtZaicSeg2}, \cite{BahtShestZ},
\cite{GMZ}, \cite{GZbook}.

For two $G$-graded $*$-algebras $A,$ $B$ their homomorphism
is called {\it graded $*$-homomorphism } if it is graded, and commutes with involution.
It happens iff it commutes with any element of $\widehat{G}$. An ideal $I \unlhd A$
of a graded algebra with involution $A$ is a graded $*$-ideal
if it is graded and invariant
under the involution. For graded algebras with involution we
consider only graded $*$-ideals and graded $*$-homomorphisms.
In this case the quotient algebra $A/I$ is also a graded algebra
with involution with the grading and involution induced from $A$.
A $G$-graded $*$-algebra is called {\it $*$-graded simple}
if it has not proper graded $*$-ideals.

We denote by $A_1 \times \dots \times A_{\rho}$ the direct product
of algebras $A_1, \dots, A_{\rho},$ and by $A_1 \oplus \dots
\oplus A_{\rho} \subseteq A$ the direct sum of subspaces $A_i$ of
an algebra $A.$ It is clear that the direct product of graded algebras
with involution is also a graded algebra with involution. Throughout the
paper we denote by $J(A)$ the Jacobson radical of $A,$ and
by $\mathrm{nd}(A)$ the degree of nilpotency of $J(A)$ if $A$
is finite dimensional. By default, all bases and dimensions of vector spaces
are considered over the base field $F.$

We always consider the lexicographical
order on the sets $\mathbb{N}_0^m,$ $m$ is a positive integer number.
Note that this order satisfies descending chain condition.

The concept of a graded identity with involution (graded $*$-identity)
is the union of concepts of a graded identity (see  \cite{BelRow},
\cite{Kem1}, \cite{GRZ1}, \cite{GZbook}) and identity with involution (see, e.g., \cite{GZbook}).
It inherits the principal features of the notion of an ordinary polynomial identity. We refer the reader
to the textbooks \cite{Dren}, \cite{DrenForm}, \cite{GZbook}, and to \cite{BelRow}, \cite{Kem1}
on questions concerning ordinary polynomial identities.

Here we study graded $*$-identities of associative $G$-graded
algebras with graded involution. We prove that a
finitely generated $G$-graded PI-algebra with involution satisfies exactly the same
graded $*$-identities as some finite dimensional graded algebra with
involution under the condition of existence of some specific basis of
a $*$-graded finite dimensional algebra (Theorem \ref{*PI-gen2}).
The required basis is defined in Lemma \ref{Pierce}.
We also give a description of $*$-graded simple finite dimensional algebras
over an algebraically closed field of characteristic zero for case
of a grading by the cyclic group $G$ of a prime order or of order $4$ (Theorem \ref{simGr-pClas}).
As a partial case we obtain \\ \\
{\bf Theorem \ref{*PI}} \ \
{\it Let $q$ be a prime integer or $q=4.$ Assume that $F$ is a field of
characteristic zero. Then for any $(\mathbb{Z}/q \mathbb{Z})$-graded
finitely generated associative PI-algebra $A$ with graded involution over $F$ there exists
a finite dimensional  over $F$ $(\mathbb{Z}/q \mathbb{Z})$-graded associative algebra $C$ with graded
involution such that the ideals of graded $*$-identities of $A$ and $C$ coincide.} \\

Finally, we suppose that the next assumption can be true in general.

\begin{conjecture}
Let $G$ be a finite abelian group, and $\widetilde{F}$ be an
algebraically closed field of characteristic zero.
Then any $G$-graded $*$-simple finite dimensional $\widetilde{F}$-algebra
possesses a basis satisfying the claims of Lemma \ref{Pierce}.
\end{conjecture}

The problem of existence of the required basis is reduced
to Assumption \ref{Class} concerning the classification of
$*$-graded simple finite dimensional algebras
over an algebraically closed field.
The confirmation of Conjecture \ref{ClFinAb} will guarantee
the existence of a basis defined in Lemma \ref{Pierce}.
In this case Theorem \ref{*PI-gen2} immediately will imply PI-representability
with respect to graded $*$-identities
of any finitely generated $G$-graded PI-algebra with graded involution
over a field $F$ of characteristic zero for any finite abelian group $G.$

\begin{conjecture}
Let $F$ be a field of characteristic zero, and $G$ a finite abelian group.
Then a proper $gi$T-ideal of graded $*$-identities  of a $G$-graded
finitely generated associative PI-algebra with graded involution over $F$ coincides
with the ideal of graded $*$-identities of some finite dimensional
over $F$ $G$-graded associative algebra with graded
involution.
\end{conjecture}

It is worth to mention that the condition for a finitely generated algebra
to be a PI-algebra in Theorems \ref{*PI-gen}, \ref{*PI-gen2}, \ref{*PI} is necessary.
Since any finite dimensional algebra is a PI-algebra (an algebra satisfying non-trivial
ordinary (non-graded) polynomial identity) then Theorems \ref{*PI-gen}, \ref{*PI-gen2}, \ref{*PI}
can be applied only to $gi$T-ideals containing some non-trivial T-ideal.
We discuss briefly in Section 1 the conditions providing that a $gi$T-ideal
contains a non-trivial T-ideal.

First, we prove Theorem \ref{*PI-gen} about PI-representability
with respect to graded $*$-identities for a field of characteristic zero which contains
a primitive root of unity of order $\mathfrak{m}=|G|.$ Afterwards, we extend this result
for any field of characteristic zero (Theorem \ref{*PI-gen2}).
In order to prove Theorem \ref{*PI-gen}
we exploit the techniques created by A.R.Kemer \cite{Kem1} for the Specht
problem solution modified for the case of graded
identities with involution. Earlier these methods also were adopted by
A.Belov e E.Aljadeff for group-graded identities \cite{AB}, and by
author for graded identities of algebras graded by a finite abelian group \cite{Svi1},
and for non-graded identities with involution \cite{Svi2}.

Here we follow the structure of the proof given in
\cite{Svi2}. The majority of constructions, properties and arguments from \cite{Svi2}
needs only slight adaptation for the graded case or even
can be directly repeated. The main definitions are given for the completeness of the text,
even if they directly repeat the non-graded versions.
The proofs are repeated only if we need to point out
some details or conditions which are peculiar in the graded case.  In all other cases
we refer the reader to the corresponding statements and arguments of \cite{Svi2}
with the appropriate comments. We can refer the reader also to
\cite{Svi1} for some technical details.

We introduce briefly in Section 1 the concept of a graded
identity with involution and the concept of the free graded algebra with involution.
In Section 2 we state the principal assumption (Assumption \ref{Class})
concerning the classification of finite dimensional $*$-graded simple algebras
over an algebraically closed field.

Section 3 is devoted to finite dimensional graded
$*$-algebras. We consider their structure, define a specific basis,
introduce the parameters $\mathrm{par}_{gi}(A)$ of a finite dimensional graded $*$-algebra $A$,
and the Kemer index $\mathrm{ind}_{gi}(\Gamma)$ of the $gi$T-ideal $\Gamma$
of graded identities with involution of a finitely generated graded PI-algebra
with involution.
We establish relations between the structural parameters
$\mathrm{par}_{gi}$ of finite dimensional graded $*$-algebras and indices $\mathrm{ind}_{gi}$
of their $gi$T-ideals. Lemmas \ref{Pierce}, \ref{ind-simple} in Section 3 are basic for the proof
and represent the principal difference with the non-graded case \cite{Svi2}. Lemma \ref{Pierce} modulo Assumption \ref{Class} substitutes Lemma 4 of the non-graded case \cite{Svi2},
and Lemma \ref{ind-simple} takes place of
Lemma 12 \cite{Svi2}. Lemmas \ref{Repr}, \ref{reduc}, \ref{Exact1}, \ref{Gammasub}
are the graded versions of Lemmas 5, 9, 14, 15 in \cite{Svi2} respectively.

Section 4 is devoted to graded identities with forms, representable algebras, and
to the technique of approximation of finitely generated graded algebras with involution
by finite dimensional graded $*$-algebras. This section almost completely repeats
the analogous section in \cite{Svi2}. Observe that the similar constructions
(the free algebra with forms, identities with forms)
can be found also in \cite{Kem2}, \cite{Rasm1}, \cite{Zubk}.
Section 5 contains the proof of the main theorems (Theorem \ref{*PI-gen}, and Theorem \ref{*PI-gen2}).
We consider algebras over a field of characteristic zero containing a primitive root of unity of order
$\mathfrak{m}=|G|$ in Theorem \ref{*PI-gen}.
Its proof is also a slight modification of the proof in the non-graded case \cite{Svi2}.
In Theorem \ref{*PI-gen2} this result is extended for case of any field of characteristic zero.

In Section 6 we consider our problem in a partial case when the group $G$ is cyclic
of a prime order $q$ or of the order $q=4.$ We give the classification
of finite dimensional $*$-graded simple algebras over an algebraically closed field
for a $(\mathbb{Z}/q \mathbb{Z})$-grading (Theorem \ref{simGr-pClas}),
and obtain the PI-representability with respect to graded $*$-identities of finitely generated
graded PI-algebras with graded involution in this case (Theorem \ref{*PI}).

Observe that in Section 1, in the definitions of free graded $*$-algebra with forms and graded $*$-identities with forms (in Section 4), and in principal Theorems \ref{*PI-gen2} (Section 5), \ref{*PI} (Section 6) we consider any field $F$ of characteristic zero.
In Section 2 and in Theorem \ref{simGr-pClas} we consider algebras over an algebraically closed field $\widetilde{F}$ of characteristic zero.
In Section 3, in the major part of Section 4, and in Theorem \ref{*PI-gen} in Section 5 we assume that $F$ contains a primitive root of unity of order $\mathfrak{m}=|G|.$

\section{Free graded algebra with involution.}

Let $F$ be a field of characteristic zero, and $G$ a finite abelian group,
$|G|=\mathfrak{m}.$
Let us consider $Y = \{ y_{i \theta} | i \in \mathbb{N},
\theta \in G \},$ \ $Z = \{ z_{i \theta}
| i \in \mathbb{N}, \theta \in G \}$ two countable
sets of pairwise different indeterminants. We denote by
$\deg_{G} y_{i \theta} = \deg_{G} z_{j \theta} =
\theta$ the $G$-degree of the variables
$Y \cup Z$ with respect to $G$-grading.
Then $Y^{\theta}=\{ y_{i \theta} | i \in \mathbb{N} \},$ \
$Z^{\theta}=\{ z_{i \theta} | i \in \mathbb{N} \}$ are homogeneous
variables of $G$-degree
$\theta \in G.$
We can define $*$-action
on monomials over $Y \cup Z$ by equalities
\begin{eqnarray} \label{freeinv}
&&w^*=(x_{i_1} \cdots x_{i_n})^*=x_{i_n}^* \cdots x_{i_1}^*
=(-1)^{\delta(w)} x_{i_n} \cdots x_{i_1}, \quad \mbox{where} \nonumber\\
&&y_{j \theta}^*=y_{j \theta}, \ \ z_{j \theta}^*=-z_{j \theta},
\quad x_j \in Y \cup Z.
\end{eqnarray}
Here the sign is determined by the parity of the number
$\delta(w)$ of variables of the set $Z$ in the monomial $w.$
Linear extension of this action is an involution on the free
associative algebra $\mathfrak{F}=F\langle Y, Z \rangle$
generated by the set $Y \cup Z.$

The algebra $\mathfrak{F}=F\langle Y, Z \rangle$ is
$G$-graded with the grading $\mathfrak{F}=\oplus_{\theta \in
G} \mathfrak{F}_{\theta}$ defined by
$\mathfrak{F}_{\theta}=\mathrm{Span}_F\{
x_{i_1} x_{i_2} \cdots x_{i_n} |
\deg_{G} x_{i_1}
\cdots \deg_{G} x_{i_n}= \theta, \ x_{j} \in Y \cup Z \}.$
It is clear that the involution (\ref{freeinv})
is graded.

The algebra $\mathfrak{F}$
is the free associative graded algebra with
involution. Its elements are called {\it graded $*$-polynomials}.
Note that the set $X^{G} = \{ x_{i \theta}=y_{i \theta}+z_{i \theta}
| i \in \mathbb{N}, \  \theta \in G \}$ generates in $\mathfrak{F}$
the $G$-graded subalgebra $F\langle X^{G} \rangle$ that is
isomorphic to the free associative $G$-graded algebra
(\cite{Svi1}).

Let $f=f(x_{1},\dots,x_{n}) \in F\langle Y, Z \rangle$ be a
non-trivial graded $*$-polynomial ($x_i \in Y \cup Z$).
We say that a graded $*$-algebra $A$
satisfies the graded $*$-identity (or graded identity with involution)
$f=0$ iff $f(a_1,\dots,a_n)=0$ for all homogeneous
$a_i \in A_{\theta_i}^{\delta_i}$ of complete degree $\deg_{\widehat{G}} a_i=\deg_{\widehat{G}} x_i=
(\delta_i,\theta_i),$ \ $\delta_i \in \{ +, - \},$ \
$\theta_i \in G$ ($i=1,\dots,n$).

Then $\mathrm{Id}^{gi}(A) \unlhd F\langle Y, Z \rangle$ is the
ideal of all graded identities with involution of $A.$ Similar to the
case of graded identities and identities with involution (\cite{Kem1},
\cite{Svi1}, \cite{Svi2}) any ideal of graded identities with
involution is two-side graded $*$-ideal of the free graded algebra with
involution $F\langle Y, Z \rangle,$ which is invariant under all its
graded $*$-endomorphisms. We call such ideals $gi$T-ideals. Also any
$gi$T-ideal $I$ of $F\langle Y, Z \rangle$ is the ideal of graded
$*$-identities of the graded algebra with involution $F\langle Y, Z \rangle/I.$
Given a set $S \subseteq F\langle Y, Z \rangle$ the
$gi$T-ideal generated by $S$ is the minimal $gi$T-ideal containing $S.$
We denote it by $giT[S]\unlhd F\langle Y, Z \rangle.$
Two $G$-graded algebras with involution $A$ and $B$ are called
{\it $gi$-equivalent } if they have the same $gi$T-ideals
of graded $*$-identities. We also say that $f=g \
(\mathrm{mod } \ \Gamma)$ for a $gi$T-ideal $\Gamma$ and
graded $*$-polynomials $f, g \in F\langle Y, Z \rangle$ if $f-g \in
\Gamma.$

We assume that a T-ideal of ordinary
polynomial identities $\Gamma_1,$ a $G$T-ideal of $G$-graded
identities $\Gamma_2$ or
a  $*$T-ideal of non-graded identities with involution
$\Gamma_3$ lies in a $gi$T-ideal $\Gamma$ if the $gi$T-ideal
$\Gamma'_i$ generated by the corresponding ideal $\Gamma_i$ lies in $\Gamma.$
Recall that a {\it PI-algebra} is an algebra
satisfying an ordinary polynomial identity.
It is clear that for a $G$-graded PI-algebra with involution $A$  the ideal of ordinary
polynomial identities $\mathrm{Id}(A),$ the ideal of graded
identities $\mathrm{Id}^{G}(A),$ and
the ideal of identities with involution
$\mathrm{Id}^{*}(A)$ lie in $\mathrm{Id}^{gi}(A).$
Note, that $A$ is a PI-algebra iff the neutral component $A_{\mathfrak{e}}$
satisfies a non-trivial $*$-identity, where $\mathfrak{e}$ is the unit element of $G$
(it follows from \cite{A3}, \cite{A4}, and \cite{BahtGR}, \cite{BergC}).
Also if $F$ contains $\sqrt[\mathfrak{m}]{1}$ then
a graded algebra with involution $A$ is PI-algebra if
it satisfies an essential $\widehat{G}$-identity (\cite{BGZ1}, see also \cite{GZbook}).

We always can assume that a generating set of a finitely generated
graded algebra with involution consists of homogeneous
elements. Given a finitely generated graded $*$-algebra $A,$ and
a finite homogeneous generating set $K$ let us denote by
$\mathrm{rkh}(K)$ the maximal number of generators of the
same complete degree in $K.$ Then $\mathrm{rkh}(A)$ is the least
$\mathrm{rkh}(K)$ for all finite homogeneous generating sets $K$ of
$A.$

We also can consider the free $G$-graded algebra with involution
$F\langle Y_{(\nu)}, Z_{(\nu)} \rangle$ of a finite rank $\nu$ generated by the sets
$Y_{(\nu)} = \{ y_{i \theta} | i=1,\dots,\nu; \theta \in G \},$ and
$Z_{(\nu)} = \{ z_{i \theta} | i=1,\dots,\nu; \theta \in G \}.$
Given a $gi$T-ideal $\Gamma \subseteq F\langle Y, Z \rangle$ and
a graded $*$-algebra $B,$ denote by $\Gamma(B)=\{ f(b_1,\dots,b_n) |
f \in \Gamma, b_i \in B \ \} \unlhd B$ the verbal ideal of $B$
corresponding to $\Gamma,$ here elements $b_i \in B$
are homogeneous of appropriate complete
degrees. Then Remark 1 \cite{Svi2} is also true in graded case.

The notions of homogeneous on polynomial degree identity, and linear identity
are analogous to the case of ordinary identities (see \cite{Dren, Kem1, GZbook}).
Similarly it is enough to consider only multilinear
graded $*$-identities in the case of characteristic zero.
Let us denote for any $\bar{n}=(n_{0 1},n_{1 1},\dots,
n_{0 \mathfrak{m}},n_{1 \mathfrak{m}}) \in \mathbb{N}_0^{2 \mathfrak{m}}$ \ ($\mathfrak{m}=|G|$)
by $P_{\bar{n}}$  the vector subspaces of $\mathfrak{F}$
formed by all multilinear polynomials depending on
$y_{1 \hat{\theta}_i},\dots,y_{n_{0 i} \hat{\theta}_i},$ \
$z_{1 \hat{\theta_i}},\dots,z_{n_{1 i} \hat{\theta_i}}$, \ $\theta_i \in G$ ($i=1,\dots,\mathfrak{m}$).
Then given a $gi$T-ideal $\Gamma$ the corresponding multilinear
parts $\Gamma_{\bar{n}}=\Gamma \cap P_{\bar{n}}$ of
$\Gamma,$ and $P_{\bar{n}}(\Gamma)=P_{\bar{n}}/\Gamma_{\bar{n}}$
of the relatively free graded algebra with involution
$\mathfrak{F}/\Gamma$ are $(FS_{n_{0 1}} \otimes \dots \otimes FS_{n_{1 \mathfrak{m}}})$-modules.
Here $S_{n_{j i}}$ acts on the corresponding set of $\widehat{G}$-homogeneous variables independently.

\begin{lemma} \label{lemma1}
Let $G$ be a finite abelian group, and $A$ a finitely generated associative
$G$-graded PI-algebra with a graded involution.
The $gi$T-ideal of graded $*$-identities of $A$ contains the ideal of graded
$*$-identities of some finite dimensional $G$-graded algebra with involution.
\end{lemma}
\noindent {\bf Proof.}
Let $\Gamma=\mathrm{Id}^{*}(A)$ be the $*$T-ideal of non-graded $*$-identities
of $A.$ Then $\Gamma$ is non-trivial and $\Gamma \subseteq \mathrm{Id}^{gi}(A).$
By \cite{Svi2} $\Gamma$ is the ideal of
$*$-identities of some finite dimensional $*$-algebra $B.$
It is clear that the $gi$T-ideal generated by $\Gamma$
coincides with the ideal of graded
$*$-identities of the finite dimensional $G$-graded algebra $B
\otimes_{F} F[G]$ with the involution induced from $B$ (and trivial on
the group algebra $F[G]$) \ $(b \otimes \theta)^*=b^* \otimes \theta$,
and the grading defined by $\deg_G b \otimes \theta = \theta$
for all $b \in B,$ $\theta \in G.$
\hfill $\Box$

\section{$*$-graded simple algebras. Assumption.}

Let $\widetilde{F}$ be an algebraically closed field of characteristic zero,
and $G$ a finite abelian group. Consider a finite dimensional
$G$-graded $\widetilde{F}$-algebra $C$ with involution. We call such an algebra
{\it $*$-graded simple} if it does not contain a non-trivial graded $*$-ideal.
It is equivalent to the condition that an algebra has no non-trivial
$\widehat{G}$-invariant ideals, where $\widehat{G}={\rm Irr}{G} \times \langle *
\rangle.$

The Jacobson radical of a finite dimensional algebra is invariant
under the action of the group of automorphisms and anti-automorphisms
$\widehat{G}.$ Then over an algebraically
closed field the radical is a $G$-graded $*$-ideal (see, e.g., \cite{GMZ}).
Particularly, any finite dimensional $*$-graded simple algebra is semisimple.

For a finite dimensional $*$-graded simple $\widetilde{F}$-algebra $C$
there exist two possibilities. Either $C$ is a $G$-graded simple algebra
with an involution compatible with the grading, or $C$ contains a proper graded
simple ideal $\mathcal{B}$. It is clear that in the second case
$C=\mathcal{B} \times \mathcal{B}^{*}.$ Hence $C$ is isomorphic to the
direct product $\mathcal{B} \times \mathcal{B}^{op}$ of a graded simple
algebra $\mathcal{B}$ and its opposite algebra $\mathcal{B}^{op}$
with the exchange involution $(a,b)^{\bar{*}}=(b,a),$ \
$a \in \mathcal{B},$ \ $b \in \mathcal{B}^{op}.$ The description
of $G$-graded simple algebras is given by the next lemma.

\begin{lemma} [Theorem 3, \cite{BahtZaicSeg}] \label{simple}
Let $\widetilde{F}$ be an algebraically closed field of characteristic zero.
Then any finite dimensional $G$-graded simple algebra $C$ over $\widetilde{F}$
is isomorphic to $M_k(\widetilde{F}^{\zeta}[H]),$ a matrix algebra over
the graded division algebra $\widetilde{F}^{\zeta}[H],$ where $H$ is a subgroup
of $G$ and $\zeta:H \times H \rightarrow \widetilde{F}^{*}$ is a 2-cocycle on $H.$
The $G$-grading on $M_k(\widetilde{F}^{\zeta}[H])$ is defined by a $k$-tuple
$(\theta_1,\dots,\theta_k) \in G^k,$ so that $\deg_G(E_{ij}
\eta_{\xi})=\theta_i^{-1} \xi \theta_j$ for any matrix unit $E_{ij}$
and any basic element $\eta_{\xi}$ of $\widetilde{F}^{\zeta}[H],$ $\xi \in H.$
\end{lemma}

Here the graded division algebra $\widetilde{F}^{\zeta}[H]=
\mathrm{Span}_F\{ \ \eta_{\xi} \ | \  \xi \in H \ \}$ is a twisted group
algebra with the product on the basic elements
$\eta_{\theta} \cdot \eta_{\xi}=\zeta(\theta,\xi) \ \eta_{\theta \xi}$
determined by a 2-cocycle $\zeta$ on a subgroup $H \leq G$ ($\theta, \xi \in H$).
It has the natural $H$-grading defined by $\deg_H \eta_{\xi}=\xi$
for any $\xi \in H.$
Note that the set  $\{ E_{ij} \eta_{\xi} | 1 \le i, j \le k, \xi \in H \}$
forms a multiplicative basis of the $G$-graded simple algebra $M_k(\widetilde{F}^{\zeta}[H]).$

\begin{definition} \label{inv-basis}
A graded involution on the $G$-graded simple algebra $M_k(\widetilde{F}^{\zeta}[H])$
is called {\it elementary} if it satisfies the condition
\begin{equation} \label{inv-bas}
(E_{ij} \eta_{\xi})^{*}=\alpha_{i,j,\xi} \  E_{i' j'} \eta_{\xi'}, \quad
1 \le i', j' \le k, \ \  \xi' \in H, \ \ \alpha_{i,j,\xi} \in \{ 1, -1\}
\end{equation}
for all $i,j=1,\dots,k,$ \ $\xi \in H.$
\end{definition}
Observe that $(i,j)=(i',j')$ in (\ref{inv-bas}) implies that $\xi=\xi',$
because the involution is graded.

Let us give the principle assumption concerning finite dimensional
$*$-graded simple algebras.

\begin{assumption} \label{Class}
We suppose that any finite dimensional $*$-graded simple algebra is isomorphic as a
graded $*$-algebra either to $G$-graded simple algebra
$\widetilde{C}^{(1)}=M_k(\widetilde{F}^{\zeta}[H])$
with an elementary involution, or to the direct product
$\widetilde{C}^{(2)}=\mathcal{B} \times \mathcal{B}^{op}$ of
a graded simple algebra $\mathcal{B}=M_k(\widetilde{F}^{\zeta}[H])$
and its opposite algebra $\mathcal{B}^{op}$ with the exchange involution $\bar{*}.$

Moreover, in any case $H$ is a subgroup of $G,$ and
$\zeta:H \times H \rightarrow \mathbb{Q}[\sqrt[\mathfrak{m}]{1}]^{*}$
is a 2-cocycle on $H$ with values in the algebraic extension of rational
numbers $\mathbb{Q}$ by a primitive root $\sqrt[\mathfrak{m}]{1}$ of one
of degree $\mathfrak{m}=|G|.$
\end{assumption}
In chapters 3, 4, 5 we consider only such $*$-graded simple algebras.

\section{Finite dimensional $*$-graded algebras.}

Let $F$ be a field of characteristic zero. Assume that $F$ contains
a primitive root $\sqrt[\mathfrak{m}]{1}$ of one of degree $\mathfrak{m}=|G|.$
Suppose that $A$ is a $G$-graded algebra with involution,
finite dimensional over the base field $F$.
Repeating the proof of Lemma 2.2
\cite{BahtZaic2} for the group $\widehat{G}={\rm Irr}{G}
\times \langle * \rangle \in Aut^{*}(A),$ and applying results of
\cite{Taft} we obtain the Wedderburn-Malcev decomposition for $G$-graded algebras
with involution.

\begin{lemma} \label{Pierce0}
Let $F$ be a field of characteristic zero containing a primitive root
$\sqrt[\mathfrak{m}]{1}$ of $1$ of degree $\mathfrak{m}=|G|.$
Any finite dimensional $G$-graded $F$-algebra with involution $A'$ is
isomorphic as a graded $*$-algebra to a $G$-graded
$F$-algebra with involution of the form
\begin{equation} \label{matrix}
A=C_1 \times \dots \times C_p \oplus J.
\end{equation}
Where the Jacobson radical $J=J(A)$ of $A$ is a graded $*$-ideal, \
$B=C_1 \times \dots \times C_p$ is a maximal semisimple graded
$*$-invariant subalgebra of $A,$ \ $C_l$ are $*$-graded simple algebras
($p \in \mathbb{N} \bigcup \{0\}$).
\end{lemma}

Given a graded $*$-algebra $B$ (not necessarily without
unit) we denote by $B^{\#}=B \oplus F \cdot 1_F$ the graded $*$-algebra
with the adjoint unit $1_F.$  We assume that $1_F$ has complete degree $(+,\mathfrak{e}).$

Similarly to \cite{Svi1, Svi2} we construct
for an algebra $A$ of the form (\ref{matrix}) the graded algebra with involution with the
free Jacobson radical. Given a graded $*$-subalgebra $\widetilde{B}
\subseteq B$ we take the free product $\widetilde{B}^{\#} *_F
F\langle Y_{(q)}, Z_{(q)} \rangle^{\#}$ of $\widetilde{B}^{\#}$ with the
free unitary graded $*$-algebra $F\langle Y_{(q)}, Z_{(q)} \rangle^{\#}$
of rank $q.$
Consider its subalgebra $\widetilde{B}(Y_{(q)}, Z_{(q)})$ generated
by the set $\widetilde{B} \cup F\langle Y_{(q)}, Z_{(q)} \rangle.$ It is
clear that $\widetilde{B}(Y_{(q)}, Z_{(q)})=\widetilde{B} \oplus
(Y_{(q)}, Z_{(q)})$
is a graded $*$-algebra, where $(Y_{(q)}, Z_{(q)})$ is the two-sided graded $*$-ideal of
$\widetilde{B}(Y_{(q)}, Z_{(q)})$ generated by the set of variables
$Y_{(q)} \cup Z_{(q)}.$

Given a $gi$T-ideal $\Gamma$ and a positive integer number $s$ consider the
quotient algebra
\begin{eqnarray} \label{FRad}
\mathcal{R}_{q,s}(\widetilde{B},\Gamma)=\widetilde{B}(Y_{(q)}, Z_{(q)})/
(\Gamma(\widetilde{B}(Y_{(q)}, Z_{(q)}))+(Y_{(q)}, Z_{(q)})^s).
\end{eqnarray}
Denote also
$\mathcal{R}_{q,s}(A)=\mathcal{R}_{q,s}(B,\mathrm{Id}^{gi}(A)).$

If $F$ is algebraically closed then using Assumption \ref{Class}
we obtain more detailed description of a finite dimensional graded
$*$-algebra.
Consider a sequence $C_1,\dots,C_p$ of $p$ $*$-graded simple algebras
of the types $\widetilde{C}^{(1)},$ $\widetilde{C}^{(2)}.$
Suppose that the algebra $C_l$ of this sequence is a
$G$-graded simple algebra with an elementary
involution (i.e. $C_l$ is of the type $\widetilde{C}^{(1)}$).
Let us denote by
\[e^{(\xi_l)}_{l,(i_l j_l)}=E_{l,i_l j_l} \eta_{\xi_l} \qquad
(1 \le i_l,j_l \le k_l, \ \  \xi_l \in H_l)\]
the basic element of $C_l=M_{k_l}(F^{\zeta_l}[H_l]).$

Consider the case when $C_l=\widetilde{C}^{(2)}=\mathcal{B} \times \mathcal{B}^{op}$
is the direct product of a $G$-graded simple algebra $\mathcal{B}=M_{k_l}(F^{\zeta_l}[H_l])$
and its opposite algebra $\mathcal{B}^{op}$ with the exchange involution
($C_l$ is of the type $\widetilde{C}^{(2)}$). Then we denote
\begin{eqnarray*}
&&e^{(\xi_l)}_{l,(i_l j_l)}=\eta_{\xi_l} (E_{l,i_l j_l},E_{l,i_l j_l})=
(E_{l,i_l j_l} \eta_{\xi_l},E_{l,i_l j_l} \eta_{\xi_l}), \quad  \mbox{ and } \\
&&\tilde{e}^{(\xi_l)}_{l,(i_l j_l)}=\eta_{\xi_l} (E_{l,i_l j_l},-E_{l,i_l j_l})=
(E_{l,i_l j_l} \eta_{\xi_l},-E_{l,i_l j_l} \eta_{\xi_l}), \quad (1 \le i_l,j_l \le k_l), \
\end{eqnarray*}
$E_{l,i_l j_l}$ are the matrix units, $\eta_{\xi_l}$ is the basic element
of $F^{\zeta_l}[H_l]$ corresponding to $\xi_l \in H_l.$
Note that in this case ($C_l=\widetilde{C}^{(2)}$) all the elements
$e^{(\xi_l)}_{l,(i_l j_l)}$
are symmetric, and $\tilde{e}^{(\xi_l)}_{l,(i_l j_l)}$
skew-symmetric with respect to involution.

For the second case let us consider also the elements
\begin{eqnarray} \label{elij}
&&e_{l,(i_l j_l,i'_l j'_l)}=\eta_{\mathfrak{e}} (E_{l,i_l j_l},E_{l,i'_l j'_l})=
(E_{l,i_l j_l} \eta_{\mathfrak{e}}, E_{l,i'_l j'_l} \eta_{\mathfrak{e}})
= \nonumber\\
&&1/2 (e^{(\mathfrak{e})}_{l,(i_l j_l)}+\tilde{e}^{(\mathfrak{e})}_{l,(i_l j_l)}
+e^{(\mathfrak{e})}_{l,(i'_l j'_l)}-
\tilde{e}^{(\mathfrak{e})}_{l,(i'_l j'_l)}), \qquad
\mbox{ and } \\
&&\tilde{e}_{l,(i_l j_l,i'_l j'_l)}=\eta_{\mathfrak{e}} (E_{l,i_l j_l},-E_{l,i'_l j'_l})=
(E_{l,i_l j_l} \eta_{\mathfrak{e}}, -E_{l,i'_l j'_l} \eta_{\mathfrak{e}})= \nonumber \\
&&1/2(e^{(\mathfrak{e})}_{l,(i_l j_l)}+\tilde{e}^{(\mathfrak{e})}_{l,(i_l j_l)}
-e^{(\mathfrak{e})}_{l,(i'_l j'_l)}+
\tilde{e}^{(\mathfrak{e})}_{l,(i'_l j'_l)}), \qquad \
(1 \le i_l,j_l,i'_l,j'_l \le k_l).  \nonumber
\end{eqnarray}
It is possible that these elements are not $G$-homogeneous,
depending on the indices $i_l, j_l, i'_l, j'_l.$

It is clear that all elements $e^{(\xi_l)}_{l,(i_l j_l)},$
$\tilde{e}^{(\xi_l)}_{l,(i_l j_l)}$ ($\xi_l \in H_l,$
$1 \le i_l,j_l \le k_l$), admitted for $C_l,$ form a $G$-homogeneous basis of $C_l.$
Moreover, their symmetric and skew-symmetric parts with respect to the involution
(eliminating proportional elements and zeros) form a $\widehat{G}$-homogeneous basis
of $C_l.$ This basis of $C_l$ will be considered as {\it canonical}.

More precisely, for simple algebras of the second type $\widetilde{C}^{(2)}$ elements
$e^{(\xi_l)}_{l,(i_l j_l)},$ $\tilde{e}^{(\xi_l)}_{l,(i_l j_l)}$ are $\widehat{G}$-homogeneous,
and their symmetric and skew-symmetric parts coincide with them, $d^{(+,\theta)}_{l i_l j_l}=e^{(\xi_l)}_{l,(i_l j_l)},$ \ $d^{(-,\theta)}_{l i_l j_l}=\tilde{e}^{(\xi_l)}_{l,(i_l j_l)},$
$\theta=\deg_G \  e^{(\xi_l)}_{l,(i_l j_l)}=\deg_G \  \tilde{e}^{(\xi_l)}_{l,(i_l j_l)}=
\theta_{i_l}^{-1} \xi_l \theta_{j_l}.$
For algebras of the type $\widetilde{C}^{(1)}$ the next alternative follows from (\ref{inv-bas}).
If $(i_l,j_l)=(i'_l,j'_l)$ then $e^{(\xi_l)}_{l,(i_l j_l)}$ is a symmetric
or skew-symmetric element. In this case $d^{(\delta,\theta)}_{l i_l j_l}=e^{(\xi_l)}_{l,(i_l j_l)}$
for $(\delta,\theta)=\deg_{\widehat{G}} e^{(\xi_l)}_{l,(i_l j_l)}.$
If $(i_l,j_l) \ne (i'_l,j'_l)$ then we have for the couple of pairs of indices
$(i_l,j_l)$ and $(i'_l,j'_l),$ and elements $\xi_l, \xi'_l \in H_l$ the equalities
$(e^{(\xi_l)}_{l,(i_l j_l)})^{*}=\alpha \  e^{(\xi'_l)}_{l,(i'_l j'_l)},$ \quad
$(e^{(\xi'_l)}_{l,(i'_l j'_l)})^{*}=\alpha \  e^{(\xi_l)}_{l,(i_l j_l)},$ where
$\alpha \in \{ 1, -1 \}.$
Then we denote the symmetric part of $e^{(\xi_l)}_{l,(i_l j_l)}$ by
$d^{(+,\theta)}_{l i_l j_l}=1/2 (e^{(\xi_l)}_{l,(i_l j_l)} \pm  e^{(\xi'_l)}_{l,(i'_l j'_l)}),$
and the skew-symmetric part by
$d^{(-,\theta)}_{l i'_l j'_l}=1/2 (e^{(\xi_l)}_{l,(i_l j_l)} \mp  e^{(\xi'_l)}_{l,(i'_l j'_l)}),$
where $\theta=\deg_G \  e^{(\xi_l)}_{l,(i_l j_l)}=\theta_{i_l}^{-1} \xi_l \theta_{j_l}=
\deg_G \  e^{(\xi'_l)}_{l,(i'_l j'_l)}=\theta_{i'_l}^{-1} \xi'_l \theta_{j'_l}.$

In any case the canonical basis of $C_l$ is formed by all non-zero elements $d^{(\delta,\theta)}_{l i_l j_l}.$
We denote by $\mathcal{I}_{l,(\delta,\theta)}=\{ (i_l,j_l) | d^{(\delta,\theta)}_{l i_l j_l} \ne 0 \ \}$
the set of all couples of indices $(i_l,j_l)$ such that the corresponding basic element
$d^{(\delta,\theta)}_{l i_l j_l}$ has the complete degree $(\delta,\theta) \in \widehat{G}.$

Thus, Lemma \ref{Pierce0} and Assumption \ref{Class} immediately imply the structure description
of a finite dimensional $*$-graded simple algebra.

\begin{lemma} \label{Pierce}
Let $F$ be an algebraically closed field of characteristic zero.
Suppose that a finite abelian group $G$ admits the classification
of finite dimensional $*$-graded simple algebras
given in Assumption \ref{Class}. Then any finite dimensional
$G$-graded $F$-algebra with involution $A$ is isomorphic
to a $G$-graded
$F$-algebra with involution $A'=C_1 \times \dots \times C_p \oplus
J.$ Where any
$*$-graded simple subalgebra $C_l$
is isomorphic to an algebra $\widetilde{C}^{(1)}$ or $\widetilde{C}^{(2)}$ of Asumption \ref{Class}
($l=1,\dots,p$).

Moreover, $A'$ can be generated as a vector space by sets of its
$\widehat{G}$-homogeneous elements $D_{(\delta,\theta)},$ \
$U_{(\delta,\theta)} \subseteq A'$
($\delta \in \{ +, -\},$ $\theta \in G$) of the form
\begin{eqnarray}
\label{basisD} &&D_{(\delta,\theta)}=\{ d^{(\delta,\theta)}_{l i_l j_l} =
\varepsilon_l d^{(\delta,\theta)}_{l i_l j_l}  \varepsilon_l \ | \  \
(i_l, j_l) \in \mathcal{I}_{l,(\delta,\theta)}; \ 1 \le l \le p \}, \\
 &&U_{(+, \theta)}=\{
(\varepsilon_{l'} r \varepsilon_{l''} + \varepsilon_{l''} r^{*}
\varepsilon_{l'})/2| \ 1 \le l' \le l'' \le p+1; \ r \in J_{\theta}
\}\nonumber
\\ \label{basisU} &&U_{(-, \theta)}=\{ (\varepsilon_{l'} r \varepsilon_{l''} -
\varepsilon_{l''} r^{*} \varepsilon_{l'})/2| \ 1 \le l' \le l''\le
p+1; \ r \in J_{\theta} \}.
\end{eqnarray}
Here $D=\bigcup_{(\delta,\theta) \in \widehat{G}} D_{(\delta,\theta)}$ is the union of the
canonical bases of $C_l$ ($l=1,\dots,p$), \  $U=\bigcup_{(\delta,\theta)
\in \widehat{G}} U_{(\delta,\theta)} \subseteq J$ is the
set of $\widehat{G}$-homogeneous radical elements.

Particularly, for any admitted element
$e^{(\xi_l)}_{l,(i_l j_l)}$ of $C_l$
($\xi_l \in G,$ $1 \le i_l,j_l \le k_l,$ $l=1,\dots,p$)
there are two possibilities. In the first case $e^{(\xi_l)}_{l,(i_l j_l)}$
is symmetric or skew-symmetric with respect to involution. Then it
coincides with the corresponding element $d^{(\delta,\theta)}_{l i_l j_l},$
\ where $\theta=\deg_G \  e^{(\xi_l)}_{l,(i_l j_l)}=\theta_{i_l}^{-1} \xi_l \theta_{j_l}.$
In the other case $e^{(\xi_l)}_{l,(i_l j_l)}$ forms a pair with the uniquely defined
element $e^{(\xi'_l)}_{l,(i'_l j'_l)}=\pm (e^{(\xi_l)}_{l,(i_l j_l)})^{*}.$
Any such pair bijectively corresponds to the pair
$\{ d^{(+,\theta)}_{l i_l j_l}, \  d^{(-,\theta)}_{l i'_l j'_l} \}$ of elements of $D.$
Where $e^{(\xi_l)}_{l,(i_l j_l)},$ \ $e^{(\xi'_l)}_{l,(i'_l j'_l)}$ are the sum, and the
difference of $d^{(+,\theta)}_{l i_l j_l},$ \  $d^{(-,\theta)}_{l i'_l j'_l}.$
And $d^{(+,\theta)}_{l i_l j_l},$ \ $d^{(-,\theta)}_{l i'_l j'_l}$ are the linear
combinations of $e^{(\xi_l)}_{l,(i_l j_l)},$ \ $e^{(\xi'_l)}_{l,(i'_l j'_l)}$
with coefficients $1/2,$ $-1/2.$ Here $\theta=\deg_G \ d^{(\delta,\theta)}_{l i_l j_l}=
\deg_G \  e^{(\xi_l)}_{l,(i_l j_l)}=\theta_{i_l}^{-1} \xi_l \theta_{j_l}.$
An element $\tilde{e}^{(\xi_l)}_{l,(i_l j_l)}$ coincides with the corresponding
$d^{(-,\theta)}_{l i_l j_l},$\  $\theta=\deg_G \ d^{(\delta,\theta)}_{l i_l j_l}=\deg_G \
\tilde{e}^{(\xi_l)}_{l,(i_l j_l)}=\theta_{i_l}^{-1} \xi_l \theta_{j_l}$ \  ($1 \le i_l,j_l \le k_l,$ $l=1,\dots,p$).

The element $\varepsilon_l=(1/\lambda_l) \sum_{i_l=1}^{k_l} e^{(\mathfrak{e})}_{l,(i_l i_l)}$
is the minimal orthogonal central idempotent
of $B'=C_1 \times \dots \times C_p,$ corresponding to the unit
element of the $l$-th $\widehat{G}$-simple component $C_l$ of the algebra
$A',$ \ $\lambda_l=\zeta_l(\mathfrak{e},\mathfrak{e}) \in
\mathbb{Q}[\sqrt[\mathfrak{m}]{1}]^{*}$ (for any $l=1,\dots,p$).

In the definition (\ref{basisU}) of the set $U=\cup_{(\delta,\theta) \in \widehat{G}}
U_{(\delta,\theta)}$ the element $r$ runs on a
$G$-homogeneous set of elements of the Jacobson
radical $J=\oplus_{l',l''=1}^{p+1} (\oplus_{\theta \in G}
\varepsilon_{l'} J_{\theta}
\varepsilon_ {l''}).$  $\varepsilon_{p+1}=1-(\varepsilon_1+ \dots +\varepsilon_p)$ is the adjoint
idempotent. Particularly, $\varepsilon_{p+1}=0$ if $A$ is a unitary algebra.
All idempotents $\varepsilon_{l}$
are $\widehat{G}$-homogeneous of degree $(+,\mathfrak{e})$ ($l=1,\dots,p+1$).
\end{lemma}

If we consider identities then the statement of Lemma \ref{Pierce}
can be extended in some sense for the case of any field
$F$ of characteristic zero containing a primitive root $\sqrt[\mathfrak{m}]{1}$.

\begin{definition}
An $F$-algebra $A$ is called representable if $A$ can be embedded
into some algebra $C$ that is finite dimensional over an extension
$\widetilde{F} \supseteq F$ of the base field $F.$
\end{definition}

\begin{lemma} \label{Repr}
Let $F$ be a field of characteristic zero containing $\sqrt[\mathfrak{m}]{1}$.
Suppose that Assumption \ref{Class} is true for any algebraically closed
extension $\widetilde{F} \supseteq F.$
Then any representable $G$-graded $F$-algebra with involution $A$ is $gi$-equivalent
to some $F$-finite dimensional $G$-graded algebra with involution $A'$ that satisfies
all the claims of Lemma \ref{Pierce}.
\end{lemma}
\noindent {\bf Proof.}
We always can assume that the extension
$\widetilde{F} \supseteq F$ is algebraically
closed. Suppose that $A$ is isomorphic to an $F$-subalgebra
$\mathcal{B}$ of a finite dimensional $\widetilde{F}$-algebra
$\widetilde{\mathcal{B}}.$ It is clear that $\mathcal{B}$ can be
considered $G$-graded with involution induced from $A.$
Consider a subalgebra $\mathcal{U}=\{ (b,b^*) \  | b \in \mathcal{B} \}$
of the $F$-algebra $\mathcal{B} \times \mathcal{B}^{op}.$
$\mathcal{U}$ is $G$-graded
with the grading $\mathcal{U}_{\theta}=\{ (b_{\theta},b^*_{\theta}) \  | b_{\theta} \in \mathcal{B}_{\theta} \},$ \ $\theta \in G,$ and has the exchange involution $(b,b^*)^{\bar{ex}}=(b^*,b),$ \ $b \in \mathcal{B}.$ Then we consider
 $\widetilde{U}=\sum_{\theta \in G} \widetilde{U}_{\theta},$
where $\widetilde{U}_{\theta}=\widetilde{F} \mathcal{U}_{\theta}
\otimes_{\widetilde{F}} \widetilde{F} \theta$ ($\theta \in G$).
$\widetilde{U}$ is an $\widetilde{F}$-subalgebra of the algebra
$(\widetilde{\mathcal{B}} \times \widetilde{\mathcal{B}}^{op})
\otimes_{\widetilde{F}} \widetilde{F}[G].$ Hence $\widetilde{U}$
is a finite dimensional $\widetilde{F}$-algebra. As an $F$-algebra
$\widetilde{U}$ is
$G$-graded with the graded involution $\star$ defined by $(\alpha (b,b^*)
\otimes \theta)^{\star}=\alpha (b^*,b) \otimes \theta,$
\ $\alpha \in \widetilde{F},$
$b \in \mathcal{B},$ $\theta \in G.$
If we consider all the algebras and graded $*$-identities over
the base field $F$ then we have
$\mathrm{Id}^{gi}(\widetilde{U})=\mathrm{Id}^{gi}(\mathcal{B})=
\mathrm{Id}^{gi}(A).$

By Lemma \ref{Pierce} the graded
$\widetilde{F}$-algebra with involution $\widetilde{U}$ has the decomposition (\ref{matrix}),
where the $*$-graded simple $\widetilde{F}$-algebras $\widetilde{C}_l$ are
of the type $\widetilde{C}^{(\mathfrak{1})}$ or $\widetilde{C}^{(\mathfrak{2})}.$
We can see from Assumption \ref{Class} that
$\widetilde{C}_{l}=\widetilde{F} C_{l}$ \ ($l=1,\dots,p$), where $C_{l}$
is the $*$-graded simple $F$-algebra generated as a vector space over the base field $F$ by the same canonical $\widehat{G}$-homogeneous basis as $\widetilde{C}_l$ over $\widetilde{F}.$
Let us take $B=C_1 \times \dots \times C_p,$ \
\ $\mathcal{R}$ is an $\widetilde{F}$-basis of $J(\widetilde{U}),$
$\Gamma=\mathrm{Id}^{gi}(A),$ $q=\dim_{\widetilde{F}}
J(\widetilde{U})=|\mathcal{R}|,$ $s=\mathrm{nd}(\widetilde{U}).$ Then the
$F$-algebra $A'=\mathcal{R}_{q,s}(B,\Gamma)$ defined by
(\ref{FRad}) is a graded algebra with
involution which is finite dimensional over $F$.
Note that $\widetilde{U}_{\theta}$ is $\widetilde{F}$-subspace
of $\widetilde{U},$ and the involution $\star$ of $\widetilde{U}$
is $\widetilde{F}$-linear. Hence, it is enough to verify multilinear graded $F$-identities
with involution of $\widetilde{U}$ only on $\widehat{G}$-homogeneous elements $b \in B,$ \
$r \in \mathcal{R}.$ It follows from the graded version of
Lemma 3 \cite{Svi2} that $A'$ satisfies all
claims of Lemma \ref{Pierce}, and $\mathrm{Id}^{gi}(A') =
\mathrm{Id}^{gi}(\widetilde{U})=\mathrm{Id}^{gi}(A).$ \hfill $\Box$

\begin{definition}\label{el-dec}
We say that an $F$-finite dimensional $G$-graded $*$-algebra $A'$ has an
{\it elementary decomposition } if it satisfies
all the claims of Lemma \ref{Pierce}.
\end{definition}

It is clear that the direct product of algebras with elementary decomposition
is the algebra with elementary decomposition. Also if $F$ is algebraically closed
and admits for the group $G$ the classification of finite dimensional $*$-graded simple algebras
given in Assumption \ref{Class} then any finite dimensional $G$-graded $F$-algebra with
involution has an elementary decomposition.

\begin{corollary} \label{eldec-fd}
Let $\widetilde{F} \supseteq F$ be an algebraically closed extension such that
Assumption \ref{Class} is true over $\widetilde{F}.$ Suppose that $F$ contains
$\sqrt[\mathfrak{m}]{1}.$ Then any finite dimensional $G$-graded $F$-algebra
with involution is $gi$-equivalent to a finite dimensional $G$-graded $F$-algebra with
involution with elementary decomposition.
\end{corollary}
\noindent {\bf Proof.}
A finite dimensional graded $F$-algebra with involution
$A$ can be naturally embedded to the graded $*$-algebra $\widetilde{A}=A \otimes_F
\widetilde{F}$ preserving graded $*$-identities. We assume here $(a \otimes \alpha)^*=
a^* \otimes \alpha,$ and $\deg_G a \otimes \alpha= \deg_G a,$
for all $a \in A,$ \ $\alpha \in \widetilde{F}.$ The algebra $\widetilde{A}$
is finite dimensional over $\widetilde{F}.$ By Lemma \ref{Repr} there
exists a finite dimensional $G$-graded $*$-algebra $A'$ with elementary decomposition
such that $\mathrm{Id}^{gi}(A') =
\mathrm{Id}^{gi}(\widetilde{A})=\mathrm{Id}^{gi}(A).$ Where all identities are considered
over the field $F.$ \hfill $\Box$

Particularly, if Assumption \ref{Class} is true for the algebraic closure
of $F$ then for graded $*$-identities of finite dimensional algebras we can consider
only algebras with elementary decomposition.

We assume further that the base field $F$ contains
a primitive root $\sqrt[\mathfrak{m}]{1}$ of one of degree $\mathfrak{m}=|G|,$ and
Assumption \ref{Class} is true for the group $G$
over any algebraically closed extension $\widetilde{F} \supseteq F.$
It means that a basis of a finite dimensional graded $F$-algebra with involution $A$ always
can be chosen in the set $D \cup U.$ Where $D$ is $\widehat{G}$-homogeneous basis of the
semisimple part $B=C_1 \times \dots \times C_p$ of $A.$
Particularly, we have $|\bigcup_{l=1}^p \mathcal{I}_{l,(\delta,\theta)}|=
|D_{(\delta,\theta)}|=\dim_F B_{\theta}^{\delta}.$

Therefore for a multilinear graded $*$-polynomial it is enough to
consider only evaluations by elements of $D \cup U$
of variables of corresponding $\widehat{G}$-degree.
Such evaluations are called {\it elementary}.
Elements of the set $D$ are called semisimple,
and elements of $U$ are radical.

Similarly to the case of group
graded identities \cite{Svi1}, and $*$-identities
\cite{Svi2} we define the numeric parameters of a finite dimensional
$G$-graded $*$-algebra, and of the ideal of graded $*$-identities of a finitely
generated $G$-graded PI-algebra with involution.
Assume that $G=\{\mathfrak{e}=\hat{\theta}_1,\hat{\theta}_2,\dots,\hat{\theta}_{\mathfrak{m}} \},$
$\mathfrak{m}=|G|,$ $\mathfrak{e}$ is the unit of $G.$ $\widehat{G}=
{\rm Irr}{G} \times \langle * \rangle \subseteq Aut^{*}(A).$

\begin{definition}
Let $A=B \oplus J$ be a finite dimensional $G$-graded $F$-algebra with graded
involution. Where $B=\sum_{(\delta, \theta) \in \widehat{G}}
B_{\theta}^{\delta}$
is a maximal semisimple graded $*$-invariant subalgebra of $A,$ and $J(A)=J$
the Jacobson radical of $A.$ We denote by
$\mathrm{dims}_{gi} A=(\dim B_{\hat{\theta_1}}^{+},
\dim B_{\hat{\theta_1}}^{-},\dots,
\dim B_{\hat{\theta}_{\mathfrak{m}}}^{+},
\dim B_{\hat{\theta}_{\mathfrak{m}}}^{-})$ the collection of dimensions of all
$\widehat{G}$-homogeneous parts of the semisimple subalgebra $B.$

Recall also that we denote by $\mathrm{nd}(A)$ the
degree of nilpotency of the radical $J.$ Then the parameter of $A$
is \ $\mathrm{par}_{gi}(A)=(\mathrm{dims}_{gi} A; \mathrm{nd}(A)).$

$\mathrm{cpar}_{gi}(A)=(\mathrm{par}_{gi}(A);\dim J(A))$
is called the complex parameter of $A.$
\end{definition}
Note that for any nonzero two-sided graded $*$-ideal $I \unlhd A$
of $A$ it holds $\mathrm{cpar}_{gi}(A/I) < \mathrm{cpar}_{gi}(A).$

Let $f=f(s_1,\dots,s_k,x_1,\dots,x_n) \in F\langle Y, Z\rangle$
be a polynomial linear on a set $S=\{ s_1,
\dots, s_k\}$ of homogeneous variables ($S \subseteq Y^{\theta}$ or
$S \subseteq Z^{\theta},$ \ $\theta \in G$). Then
$f$ is alternating on $S,$ if
$f(s_{\sigma(1)},\dots,s_{\sigma(k) },x_1,\dots,x_n)=(-1)^{\sigma}
f(s_1,\dots,s_k,x_1,\dots,x_n)$ holds for any permutation $\sigma
\in \mathrm{S}_k.$

It is clear that a polynomial
$f(s_1,\dots,s_k,x_1,\dots,x_n)$ is alternating on the set $S=\{ s_1, \dots, s_k\}$ if and only if
$$f(s_1,\dots,s_k,x_1,\dots,x_n)=
\mathcal{A}_{S}(g)=\sum_{\sigma \in \mathrm{S}_k} (-1)^{\sigma}
g(s_{\sigma(1)},\dots,s_{\sigma(k)},x_1,\dots,x_n)$$ for some graded
$*$-polynomial $g(s_1,\dots,s_k,x_1,\dots,x_n)$
which is linear on the set $S.$
The mapping $\mathcal{A}_{S}$ is a graded linear transformation
and is called alternator. The properties of alternating graded
$*$-polynomials are similar to the case of ordinary polynomials (see, e.g.
\cite{Dren}, \cite{Kem1}, \cite{GZbook}).
Note also that an alternator commutes with involution.

Given a $2 \mathfrak{m}$-tuple
$\overline{t}=(t_1,\dots,t_{2 \mathfrak{m}}) \in
\mathbb{N}_0^{2 \mathfrak{m}}$ we
say that a graded $*$-polynomial $f \in F \langle Y, Z \rangle$ has the
collection of $\overline{t}$-alternating homogeneous variables
($f$ is $\overline{t}$-alternating) if
$f(Y_1,Z_1,\dots,Y_{\mathfrak{m}},
Z_{\mathfrak{m}},X)$ is linear on
$\bigcup_{j=1}^{\mathfrak{m}} (Y_{j} \cup Z_{j}),$
and $f$ is alternating on each of the sets
$Y_{j} \subseteq Y^{\hat{\theta}_j},$ \
$Z_{j} \subseteq Z^{\hat{\theta}_j},$
where $|Y_{j}|=t_{2 j-1},$ \ $|Z_{j}|=t_{2 j},$ \
$j=1,\dots,\mathfrak{m}.$

Recall that we order $2 \mathfrak{m}$-tuples lexicographically.
The definitions of the type of a multihomogeneous on degrees polynomial $f \in F \langle Y, Z \rangle,$
and the Kemer index of $gi$T-ideal of a finitely generated graded PI-algebra with graded involution
repeat the corresponding definitions for the case of graded polynomials
\cite{Svi1}. We will consider them for the completeness of the text.

\begin{definition}
Given a $2 \mathfrak{m}$-tuple $\overline{t}=(t_1,\dots,t_{2 \mathfrak{m}})
\in \mathbb{N}_0^{2 \mathfrak{m}}$
consider some (possibly different) collections  $\tau_1, \dots, \tau_s \in
\mathbb{N}_0^{2 \mathfrak{m}}$ satisfying the conditions
$\tau_j > \overline{t}$ for any $j=1, \dots, s.$
Let $f \in F\langle Y, Z \rangle$ be a multihomogeneous graded $*$-polynomial. Suppose that
$f=f(S_1,\dots,S_{s+\mu};X)$ is simultaneously $\tau_j$-alternating on
$S_j=\cup_{\theta \in G} ( Y_{j}^{\theta} \cup Z_{j}^{\theta})$
for any $j=1,\dots,s,$ and $\overline{t}$-alternating on any
$S_j=\cup_{\theta \in G} (Y_{j}^{\theta} \cup Z_{j}^{\theta}),$
$j=s+1,\dots,s+\mu.$ All the collections $S_j$ are pairwise disjoint.
Then we say that $f$ has the type
$(\overline{t};s;\mu).$ Here $|Y_{j}^{\hat{\theta}_i}|=\tau_{j, 2 i-1},$
$|Z_{j}^{\hat{\theta}_i}|=\tau_{j, 2 i}$
for any $j=1,\dots,s$ or $|Y_{j}^{\hat{\theta}_i}|=t_{2 i-1},$
$|Z_{j}^{\hat{\theta}_i}|=t_{2 i}$
\ for all $j=s+1,\dots,s+\mu$ ($i=1,\dots,\mathfrak{m}$).
\end{definition}
Note that multihomogeneous polynomials $f$ and $f^*$
always have the same type.

\begin{definition} \label{defbeta}
Given a $gi$T-ideal $\Gamma \unlhd F\langle Y, Z \rangle$ the
parameter $\beta(\Gamma)=\overline{t}$ is the greatest
lexicographic $2 \mathfrak{m}$-tuple $\overline{t} \in
\mathbb{N}_0^{2 \mathfrak{m}}$
such that for any $\mu \in \mathbb{N}$ there exists a graded
$*$-polynomial $f \notin \Gamma$ of the type $(\overline{t};0;\mu).$
\end{definition}

\begin{definition} \label{defgamma}
Given a nonnegative integer $\mu$ let $\gamma(\Gamma;\mu)=s \in
\mathbb{N}$ be the smallest natural $s$ such that any
graded $*$-polynomial of the type $(\beta(\Gamma);s;\mu)$ belongs to
$\Gamma.$

$\gamma(\Gamma;\mu)$ is a positive non-increasing function on
$\mu.$ Let us denote the limit of this function
$\gamma(\Gamma)=\lim \limits_{\mu \to \infty} \gamma(\Gamma;\mu)
\in \mathbb{N}.$
\end{definition}

\begin{definition}
The $(2 \mathfrak{m}+1)$-tuple
$\mathrm{ind}_{gi}(\Gamma)=(\beta(\Gamma); \gamma(\Gamma))$
is called by Kemer index of a $gi$T-ideal $\Gamma.$
\end{definition}

A finitely generated PI-algebra satisfies a non-graded Capelli identity
\cite{Kem6}. Similarly to the case of $G$T-ideals \cite{Svi1} and $*$T-ideals \cite{Svi2}
the Kemer index is well defined for the $gi$T-ideal of graded $*$-identities
of a finitely generated $G$-graded PI-algebra with involution.
We denote
$\mathrm{ind}_{gi}(A)=\mathrm{ind}_{gi}(\mathrm{Id}^{gi}(A))$ for a
finitely generated graded PI-algebra $A$ with involution.
$A$ is nilpotent of degree $s$ if and only if
$\mathrm{ind}_{gi}(A)=\mathrm{par}_{gi}(A)=(0,\dots,0;s).$

We obtain automatically the notion of $\mu$-boundary polynomials for
a $gi$T-ideal.
\begin{definition} [Definition 7, \cite{Svi2}] \label{Kemerpolyn}
Given a nonnegative integer $\mu$ a nontrivial multihomogeneous
polynomial $f \in F\langle Y, Z \rangle$ is called $\mu$-boundary polynomial
for a $gi$T-ideal $\Gamma$ if $f \notin \Gamma,$ and $f$
has the type $(\beta(\Gamma);\gamma(\Gamma)-1;\mu).$

Denote by $S_\mu(\Gamma)$ the set of all $\mu$-boundary
polynomials for $\Gamma$. Then
$S_\mu(A)=S_\mu(\mathrm{Id}^{gi}(A)),$ \
$K_\mu(\Gamma)=giT[S_\mu(\Gamma)],$ \ $K_{\mu,
A}=giT[S_\mu(A)].$
\end{definition}

The set  $S_\mu(\Gamma)$
of all $\mu$-boundary polynomials of a $gi$T-ideal $\Gamma,$
and the Kemer index satisfy the same basic properties as in the case
of $G$T-ideals and $*$T-ideals (see Lemmas 4-10 \cite{Svi1}, and Lemmas 6-8 \cite{Svi2}).
Observe that these properties do not depend on the type of identities.
They are completely determined by Definitions \ref{defbeta}-\ref{Kemerpolyn}
(see the arguments in \cite{Svi1}).

We can consider also the graded version of a $*$PI-reduced algebra.
We call it $gi$-reduced algebra.
\begin{definition}
A finite dimensional $G$-graded $*$-algebra $A$ with elementary decomposition is
$gi$-reduced if there do not exist finite dimensional
$G$-graded $*$-algebras with elementary decomposition
$A_1,\dots,A_\varrho$ ($\varrho \in \mathbb{N}$) such that
$\bigcap\limits_{i=1}\limits^{\varrho} \mathrm{Id}^{gi}(A_i) =
\mathrm{Id}^{gi}(A),$ and $\mathrm{cpar}_{gi}(A_i) < \mathrm{cpar}_{gi}(A)$
for all $i=1,\dots,\varrho.$
\end{definition}

The next graded modification of Lemma 9 \cite{Svi2} holds.

\begin{lemma} \label{reduc}
Given a $gi$-reduced algebra $A$ with the Wedderburn-Malcev
decomposition (\ref{matrix}) $A=(C_1 \times \cdots \times C_p)
\oplus J,$ we have $C_{\sigma(1)} J C_{\sigma(2)} J \cdots J
C_{\sigma(p)} \ne 0$ for some $\sigma \in \mathrm{S}_p.$
\end{lemma}
\noindent {\bf Proof.}
Suppose that $C_{\sigma(1)} J C_{\sigma(2)} J \cdots J C_{\sigma(p)} = 0$
for any $\sigma \in \mathrm{Sym}_p.$ Consider
$G$-graded $*$-algebras with elementary decomposition
$A_i=(\prod \limits_{\mathop{1 \le j \le
p}\limits_{\scriptstyle j \ne i}} C_j ) \oplus J(A)$ ($i=1,\dots, p$).
We have for them $\mathrm{Id}^{gi}(A)=\bigcap
\limits_{i=1}\limits^{p} \mathrm{Id}^{gi}(A_i),$ and $\mathrm{dims}_{gi}
A_i < \mathrm{dims}_{gi} A$ for any $i=1,\dots,p.$
This contradicts to the definition of $gi$-reduced algebra.
\hfill $\Box$

Particularly, we have $\mathrm{nd}(A) \ge p$ for a $gi$-reduced algebra $A.$
Then Corollary \ref{eldec-fd} along with the properties of $\mu$-boundary polynomials implies
also the graded versions of Lemmas 10, 11 \cite{Svi2} (see the proofs in \cite{Svi1}).

The Kemer index and parameters of $gi$-reduced algebras are related in a similar way as in case of non-graded involution. It is the crucial point of our proof.

\begin{lemma} \label{ind-simple}
Given a $gi$-reduced algebra $A$ we have
$\beta(A)=\mathrm{dims}_{gi} A.$
Any $*$-graded simple finite dimensional
algebra $C$ with elementary decomposition is $gi$-reduced, and
$\mathrm{ind}_{gi}(C)=\mathrm{par}_{gi}(C)=(t_1,\dots,t_{2 \mathfrak{m}};1).$
\end{lemma}
\noindent {\bf Proof.}
If $A$ is nilpotent then the assertion of Lemma is trivial.
Suppose that $A$ is a non-nilpotent $gi$-reduced algebra.
By the graded version of Lemma 6 \cite{Svi2} we have $\beta(A) \le \mathrm{dims}_{gi} A.$
Thus it is enough to find a graded
$*$-polynomial of the type $(\mathrm{dims}_{gi} A;0;\hat{s})$
which is not a graded identity with involution of $A$ for any $\hat{s} \in \mathbb{N}.$

Consider the elementary decomposition (\ref{matrix}) of $A.$
Similarly to Lemma 12 \cite{Svi2}
for any $*$-graded simple component $C_l$
($l=1,\dots,p$) we take $\hat{s}$ sets of
distinct $\widehat{G}$-homogeneous variables
$Y^{(\delta,\theta)}_{l, m}=\{ y^{(\delta,\theta)}_{l, (i_l j_l), m} | (i_l,j_l)
\in \mathcal{I}_{l,(\delta,\theta)} \}$ corresponding to the canonical basic elements
$d^{(\delta,\theta)}_{l i_l j_l} \in D_{(\delta,\theta)}$ (see (\ref{basisD})).
Here $Y^{(+,\theta)}_{l, m} \subseteq Y^{\theta},$
and $Y^{(-,\theta)}_{l, m} \subseteq Z^{\theta},$ \ $\delta \in \{ +, -\},$
$\theta \in G,$ \ $m=1,\dots,\hat{s}.$

Suppose that $C_l=\widetilde{C}^{(\mathfrak{q})}$
with $\mathfrak{q}=1,2.$
Then consider the polynomial $w_{l, m}(Y_{l, m},X_{l})$ which is the product
of all variables of the set $Y_{l, m}=\cup_{(\delta,\theta) \in \widehat{G}}
Y^{(\delta,\theta)}_{l, m}$ connected by $x_{l,(i j)}$ if
$\mathfrak{q} = 1$ or by $x_{l, (i j,i' j')}$ if
$\mathfrak{q} = 2.$
Here we take $x_{l,(i j)}=\pi_{1} \ \tilde{y}_{l,(i j)}  +
\pi_{2} \ \tilde{z}_{l,(i j)},$ and
$x_{l, (i j,i' j')}=1/2(\pi_{3} \  \tilde{y}_{l,(i j)} +
\pi_{4} \ \tilde{z}_{l,(i j)} + \pi_{5} \  \tilde{y}_{l,(i' j')} +
\pi_{6} \  \tilde{z}_{l,(i' j')})$ for $1 \le i, j, i', j' \le k_l.$
Where $\pi_{s} \in \{ 0, 1, -1\},$ and\ $\tilde{y}_{l,(i j)} \in Y,$ \
$\tilde{z}_{l,(i j)} \in Z$ with $\deg_{G} \tilde{y}_{l,(i j)} =
\deg_{G} \tilde{z}_{l,(i j)} = \deg_{G}
e^{(\mathfrak{e})}_{l,(i j)}.$ We also say here that
$x_{l, (j_1 i_2)},$ and  $x_{l, (j_1 i_2,j_2 i_1)}$ connect the
variable $y^{(\delta_1,\theta_1)}_{l, (i_1 j_1), m}$ with
$y^{(\delta_2,\theta_2)}_{l, (i_2 j_2), m},$ and $y^{(\delta_1,\theta_1)}_{l, (i_1 j_1), m} x_{l, (j_1 i_2)} y^{(\delta_2,\theta_2)}_{l, (i_2 j_2), m},$ \  $y^{(\delta_1,\theta_1)}_{l, (i_1 j_1), m} x_{l, (j_1 i_2,j_2 i_1)} y^{(\delta_2,\theta_2)}_{l, (i_2 j_2), m}$ are the corresponding
connected product of these variables. Then we denote $X_l = \{ \tilde{y}_{l,(i j)},
\tilde{z}_{l,(i j)}| 1 \le i, j \le k_l \}.$

By Lemma \ref{reduc} we can assume without lost of generality that $A$ contains an element
\begin{equation} \label{ErE}
a=e^{(\mathfrak{e})}_{1,(s_1 s_1)} r_1 e^{(\mathfrak{e})}_{2,(s_2 s_2)} \cdots e^{(\mathfrak{e})}_{p-1,(s_{p-1}
s_{p-1}) } r_{p-1} e^{(\mathfrak{e})}_{p,(s_{p} s_{p}) } \ne 0,
\end{equation}
where $r_l \in U$ are some $\widehat{G}$-homogeneous basic radical elements
chosen in the set $U$ (Lemma \ref{Pierce}).

Then in the case $\mathfrak{q} = 1$ let us consider the graded $*$-polynomial
$W_{l,s_l}=x_{l, (s_l t_l)}\cdot \bigl(
\prod_{m=1}^{\hat{s}}  (x_{l, (t_l 1)} \cdot w_{l, m}) \bigr)
 \cdot x'_{l, (t_l s_l)}.$
And for $\mathfrak{q}=2$ let us take
$W_{l,s_l}=x_{l, (s_l t_l,t'_l s_l)}\cdot \bigl(
\prod_{m=1}^{\hat{s}}  (x_{l, (t_l 1,1 t'_l)} \cdot w_{l, m}) \bigr)
 \cdot x'_{l, (t_l s_l,s_l t'_l)}.$ Where $s_l$ are given by
(\ref{ErE}), and $t_l,$ $t'_l$ are chosen to connect the word $w_{l, m}$
with $w_{l, m+1}.$ The variables
$x'_{l, (i j)},$ and $x'_{l, (i j,i' j')}$ are defined
as linear combinations of
$\widehat{G}$-homogeneous variables of the set $X'_l = \{ \tilde{y}'_{l,(i j)},
\tilde{z}'_{l,(i j)}| 1 \le i, j \le  k_l \}$ similarly to $x_{l, (i j)},$ and
$x_{l, (i j,i' j')}.$

Denote $Y^{(\delta,\theta)}_{(m)}= \bigcup_{l=1}^{p}
Y^{(\delta,\theta)}_{l,m}$ (${(\delta,\theta)} \in \widehat{G}$).
The polynomial
\begin{eqnarray} \label{polyn2}
f=\bigl( \prod_{m=1}^{\hat{s}} \prod_{(\delta,\theta) \in \widehat{G}}
\mathcal{A}_{Y^{(\delta,\theta)}_{(m)}} \bigr) \  W_{1,s_1} \  \hat{x}_1 \  W_{2,s_2} \  \hat{x}_2 \cdots \hat{x}_{p-1} \   W_{p,s_p}
\end{eqnarray}
is $(\mathrm{dims}_{gi} A)$-alternating on any set $Y_{(m)}$ for all $m=1,\dots,\hat{s}.$
Here $\hat{x}_q \in Y^{\theta}$ if
$r_{q} \in U_{(+,\theta)},$ and $\hat{x}_q \in Z^{\theta}$ if
$r_{q} \in U_{(-,\theta)}$ in (\ref{ErE}), $q=1,\dots,p-1.$

Then the evaluation
\begin{eqnarray} \label{s-sim}
&&y^{(\delta,\theta)}_{l,(i_l j_l), m} = d^{(\delta,\theta)}_{l i_l j_l},  \qquad \
(i_l,j_l) \in \mathcal{I}_{l,(\delta,\theta)}; \nonumber \\
&&\hat{x}_{q}=r_q; \nonumber
\\&&l=1,\dots,p; \ \ q=1,\dots, p-1; \ \ m=1,\dots,\hat{s}; \ \
(\delta,\theta) \in \widehat{G};
\end{eqnarray}
of the polynomial $f$ is equal to $\alpha a \ne 0,$
where the element $a \in A$ is defined in (\ref{ErE}),
and the non-zero coefficient $\alpha \in F$ is the product
of the corresponding values of the 2-cocycles $\zeta_l$ divided by $2^{c_l},$ \ $c_l \in \mathbb{N}$
(see the proof of Lemma 12 \cite{Svi2}, and Lemmas 11, 15 \cite{Svi1}).
Here the elementary substitution of variables of the sets $X_l,$
$X'_{l},$ and the coefficients $\pi_{s} \in \{0, 1, -1 \}$
for any collection $(l,(i,j))$
are chosen to guarantee $x_{l, (i_l j_l)}=
e^{(\mathfrak{e})}_{l,(i_l j_l)},$ \
$x_{l, (i_l j_l,i'_l j'_l)}=e_{l,(i_l j_l,i'_l j'_l)}$
(see Lemma \ref{Pierce}, and (\ref{elij}));
and \ $x'_{l, (i_l j_l)}=e^{(\theta_l)}_{l,(i_l j_l)},$
\  $x'_{l, (i_l j_l,i'_l j'_l)}=\eta_{\theta_l} (E_{l i_l j_l},(-1)^{k_l^2 |H_l|} E_{l i'_l j'_l})$
for the suitable element $\theta_l \in H_l$
($1 \le i_l, j_l, i'_l, j'_l \le k_l$, $1 \le l \le p$).

Therefore we have $f \notin
\mathrm{Id}^{gi}(A).$ Hence at
least one multihomogeneous component $\tilde{f}$ of $f$
is not a graded $*$-identity of $A$ also, and it is
$(\mathrm{dims}_{gi} A)$-alternating on any set
$Y_{(m)},$ $m=1,\dots,\hat{s}.$
Thus $\tilde{f}$ is the required polynomial.

Notice that the condition of $gi$-reducibility of $A$ is necessary
only to find in $A$ a non-zero element (\ref{ErE}).
If $A$ is a $*$-graded simple algebra with elementary decomposition ($p=1$) then
instead of $a$ in (\ref{ErE}) we take
$e^{(\mathfrak{e})}_{1,(1 1)},$ and
will also obtain $\beta(A)=\mathrm{dims}_{gi} A.$
Since $\mathrm{ind}_{gi}(A) \le \mathrm{par}_{gi}(A)=(\beta(A);1),$
and $\gamma(A)>0$ then
$\mathrm{ind}_{gi}(A) = \mathrm{par}_{gi}(A).$ By the graded version of Lemma 6 \cite{Svi2} the
conditions $\dim J(A)=0,$ and $\mathrm{ind}_{gi}(A) =
\mathrm{par}_{gi}(A)$ imply that $A$ is $gi$-reduced.
 \hfill $\Box$

Assume that a finite dimensional graded algebra with involution $A$
has an elementary decomposition.
Similarly to \cite{Svi1}, \cite{Svi2} we define special types
of evaluations of variables of a graded $*$-polynomial in $A$: elementary complete
and elementary thin evaluations. We also define the notion of an exact polynomial
related with these evaluations.

\begin{definition}
An elementary evaluation $(a_1,\dots,a_n)$ of $\widehat{G}$-homogeneous elements of $A$
(namely, $a_i \in D \cup U \subseteq A$
(\ref{basisD}), (\ref{basisU})) is called incomplete if there
exists $j=1,\dots,p$ such that $$\{a_1,\dots,a_n\} \cap \bigl(C_j
\oplus_{l=1}^{p+1} (\varepsilon_j J \varepsilon_l + \varepsilon_l
J \varepsilon_j) \bigr) = \emptyset.$$ Otherwise the evaluation
$(a_1,\dots,a_n)$ is called complete.
\end{definition}

\begin{definition}
An elementary evaluation $(a_1,\dots,a_n) \in A^n$ is
called thin if it contains strictly less than $\mathrm{nd}(A)-1$
radical elements (not necessarily distinct).
\end{definition}

\begin{definition}
We say that a multilinear graded $*$-polynomial
$f(x_1,\dots,x_n) \in F\langle Y, Z \rangle$
is exact for a finite dimensional graded $*$-algebra $A$ with elementary decomposition
if $f(a_1,\dots,a_n)=0$ holds in $A$ for any thin or incomplete
evaluation $(a_1,\dots,a_n) \in A^n.$
\end{definition}

By Lemma \ref{Pierce} we have that
$\dim_F (C_{l})_{\mathfrak{e}}^{+}
> 0$ for any $l=1,\dots,p$. Hence the arguments similar to original ones
prove the graded version of Lemma 13 \cite{Svi2}.
The next two statements are also true for graded $*$-identities.

\begin{lemma} \label{Exact1}
Any nonzero $gi$-reduced algebra $A$ has an exact polynomial,
that is not a graded $*$-identity of $A.$
\end{lemma}
\noindent {\bf Proof.} If $A$ is nilpotent then the assertion
follows from the graded version of Lemma 13 \cite{Svi2}. Suppose that $A$ is a non-nilpotent
graded $*$-algebra satisfying the claims of Lemma \ref{Pierce}.
Consider its $*$-invariant graded subalgebras $A_i=(\prod
\limits_{\mathop{1 \le j \le p}\limits_{\scriptstyle j \ne i}} C_j
) \oplus J(A)$ for all $i=1,\dots, p.$ Take $q=\dim_F J(A),$ \
$s=\mathrm{nd}(A)-1.$
By the graded version of Lemma 3 \cite{Svi2} the graded $*$-algebra
$\mathcal{R}_{q,s}(A)$ has an elementary decomposition, and
$\mathrm{Id}^{gi}(A) \subseteq \mathrm{Id}^{gi}(\widetilde{A}),$
where  $\widetilde{A}=A_1 \times \dots
\times A_p \times \mathcal{R}_{q,s}(A).$
Consider any map of the set
$Y_{(q)} \cup Z_{(q)}$
onto a $\widehat{G}$-homogeneous basis of $J(A)$ of the form (\ref{basisU}) which
preserves $\widehat{G}$-degrees of variables. Such map can be extended to a surjective graded
$*$-homomorphism $\varphi: B(Y_{(q)},Z_{(q)}) \rightarrow A,$ assuming
$\varphi(b)=b$ for any $b \in B.$

Therefore similarly to Lemma 14 \cite{Svi2} any multilinear graded $*$-identity
of the algebra $\widetilde{A}$ is exact for $A.$
It is clear that all the algebras $A_i,$ and the algebra
$\mathcal{R}_{q,s}(A)$ have the complex parameter less than $A.$
Since $A$ is $gi$-reduced then we have $\mathrm{Id}^{gi}(A)
\subsetneqq \mathrm{Id}^{gi}(\widetilde{A})$.
Any multilinear polynomial $f$ such that $f \in \mathrm{Id}^{gi}(\widetilde{A}),$
and $f \notin \mathrm{Id}^{gi}(A)$ satisfies the assertion of the lemma. \hfill $\Box$

\begin{lemma} \label{Gammasub}
Let $A$ be a finite dimensional graded $*$-algebra with an elementary decomposition,
$h$ an exact polynomial for $A,$ and $\bar{a} \in A^n$ a complete
evaluation of $h$ containing exactly $\tilde{s}=\mathrm{nd}(A)-1$
radical elements. Then for any $\mu
\in \mathbb{N}_0$ there exist a graded $*$-polynomial
$h_\mu \in giT[h],$ and
an elementary evaluation $\bar{u}$ of $h_\mu$ in $A$  such that:
\begin{enumerate}
 \item $h_\mu(\mathcal{Z}_1,\dots,\mathcal{Z}_{\tilde{s}+\mu},
\mathcal{X})$ is
$\tau_j$-alternating on any set $\mathcal{Z}_j$ with $\tau_j >
\beta=\mathrm{dims}_{gi} A$ for all $j=1,\dots,\tilde{s}$, and
is $\beta$-alternating on any $\mathcal{Z}_j$ for
$j=\tilde{s}+1,\dots,\tilde{s}+\mu$ (all the sets
$\mathcal{Z}_j,$ $\mathcal{X} \subseteq (Y \cup Z)$ are disjoint),
 \item $h_\mu(\bar{u})= h(\bar{a}),$
 \item all the variables from $\mathcal{X}$ are
replaced by semisimple elements.
\end{enumerate}
\end{lemma}
\noindent {\bf Proof.}
Consider the decomposition (\ref{matrix}) of $A.$
Take any $l=1,\dots,p.$ Let $W_{l,s_l}(\widetilde{Y}_l,\widetilde{X}_l)$ be defined as
in Lemma \ref{ind-simple} for $\hat{s}=\tilde{s}+\mu.$
Here $\widetilde{Y}_l = \cup_{m=1}^{\tilde{s}+\mu} Y_{l,m},$ \
$\widetilde{X}_l=X_l \cup X'_l.$ Suppose that the evaluation (\ref{s-sim}) of the polynomial
$W_{l,s_l}$ is equal to $\alpha_{l,s_l} e^\mathfrak{e}_{l, (s_l,s_l)}$
(Lemma \ref{ind-simple}, see also Lemma 12 \cite{Svi2},
and Lemmas 11, 15 \cite{Svi1}), where $\alpha_{l,s_l} \in F$ is the non-zero coefficient.
Consider the polynomial
$\tilde{f}_l(\widetilde{Y}_l,\widetilde{X}_l)=(1/\lambda_l) \sum_{s_l=1}^{k_l}
{1 / \alpha_{l,s_l}} W_{l,s_l},$
where $\lambda_l=\zeta_l(\mathfrak{e},\mathfrak{e})$ (see Lemma \ref{Pierce}).

Notice that the polynomial $\tilde{f}_l$ is not
necessary $G$-homogeneous due to terms $x_{l, (i j,i' j')}$ that can be non-homogeneous.
Denote by $\bar{f}_l$ the $\mathfrak{e}$-component of
$\tilde{f}_l$ in $G$-grading, and by $f^{\mathfrak{e}}_l=(\bar{f}_l+
\bar{f}_l^{*})/2$ its symmetric part.

From the proof of Lemma \ref{ind-simple} it is clear that the evaluation (\ref{s-sim}) of the polynomial
$\tilde{f}'_l=(\prod_{m=1}^{\tilde{s}+\mu} \prod_{(\delta,\theta) \in \widehat{G}}
\mathcal{A}_{Y^{(\delta,\theta)}_{l,m}}) \tilde{f}_l$ is equal to $\varepsilon_l.$
Since $\varepsilon_l$ is a $\widehat{G}$-homogeneous element of degree
$(+,\mathfrak{e})$ then the result of
this evaluation of the polynomial $f'_l=(\prod_{m=1}^{\tilde{s}+\mu} \prod_{(\delta,\theta) \in \widehat{G}}
\mathcal{A}_{Y^{(\delta,\theta)}_{l,m}}) f^{\mathfrak{e}}_l$
is the same. Recall that any alternator is graded and commutes with involution.

Assume that $\zeta_1, \dots, \zeta_{\tilde{s}}$ are the variables of $h$ evaluated by radical elements of $\bar{a}.$  Let us denote
$\mathcal{Z}^{(\delta,\theta)}_m=\bigcup_{l=1}^p (Y^{(\delta,\theta)}_{l, m}) \cup \{ \zeta_m \}$
if $m=1,\dots,\tilde{s},$ and $\deg_{\widehat{G}} \zeta_m = (\delta,\theta)$ or $\mathcal{Z}^{(\delta,\theta)}_m=\bigcup_{l=1}^p (Y^{(\delta,\theta)}_{l, m})$
otherwise. Then $\mathcal{Z}_m=\bigcup_{(\delta,\theta) \in \widehat{G}} \mathcal{Z}^{(\delta,\theta)}_m.$
Let us denote by $\hat{f}_1$ the polynomial $\frac{1}{2} f^{\mathfrak{e}}_1$ in the case $p=1,$ and $h(\bar{a}) \in \varepsilon_1 A \varepsilon_1.$ The polynomials $\hat{f}_l=f^{\mathfrak{e}}_l$ ($l=1,\dots,p$) must be taken in all other cases.
We obtain the proof of our lemma in graded case if we replace the polynomials $f_l$ by $\hat{f}_l$
in the proof of Lemma 15 \cite{Svi2} and apply to the polynomial $h'$ the product of the alternators acting on $\mathcal{Z}^{(\delta,\theta)}_m$
(for all $(\delta,\theta) \in \widehat{G},$  \  $m=1,\dots,\tilde{s}+\mu$).
Note that all another elements considered in Lemma 15 \cite{Svi2} are $G$-homogeneous.
Remark also that we obtain the evaluation $\bar{u}$ replacing the variables of the polynomials $\hat{f}_l$ as in Lemma \ref{ind-simple} (see (\ref{s-sim})), and the variables $\zeta_m,$ \
$\tilde{x}_{n'}$ by the corresponding elements $a_s$ as in Lemma 15 \cite{Svi2} (see (13)).
\hfill $\Box$

Similarly to the case of non-graded $*$-polynomials \cite{Svi2}, and
polynomials graded by an abelian group \cite{Svi1} Lemmas \ref{ind-simple}, \ref{Exact1}, \ref{Gammasub}
imply the graded versions of Lemmas 16-19 \cite{Svi2}.

\section{Representable algebras.}

Consider a graded version of $*$-identities with forms introduced in
\cite{Svi2}.
Let $F$ be a field, and $R$ a commutative associative $F$-algebra.
Suppose that a $G$-graded $F$-algebra $A$
with involution has a structure of $R$-algebra satisfying
$R A_\theta \subseteq A_\theta,$ \ $\forall \theta \in G,$
and the involution of $A$ is $R$-linear, i.e.
$\deg_{G} r a = \deg_{G} a,$ \
$r a=a r,$ \ $(ra)^*=r a^*$ for all $r \in R,$ \ $a \in A.$
Particularly, it happens if $R=F$ or if $R \subseteq Z(A) \cap
A_{\mathfrak{e}}^{+},$
where $Z(A)$ is the center of $A.$

\begin{definition} [Definition 13 \cite{Svi2}]
Let $A$ be an $R$-algebra with involution.
Any $R$-multilinear mapping $\mathfrak{f}:A^{n} \rightarrow R$ is
called $n$-linear form on $A.$
\end{definition}

\begin{definition}
Suppose that $A,$ $B$ are $F$-algebras with an $n$-linear form  $\mathfrak{f}.$
A homomorphism of $F$-algebras $\varphi:A \rightarrow B$ preserves
the form $\mathfrak{f}$ if
\[ \varphi(a_0 \mathfrak{f}(a_1,\dots,a_n)) =
\varphi(a_0) \mathfrak{f}(\varphi(a_1),\dots,\varphi(a_n)),  \quad \forall a_i \in A.\]
\end{definition}

Let us consider the free $*$-algebra with forms
$FS\langle Y, Z \rangle = F\langle Y, Z \rangle \otimes_F \mathcal{S}$
defined for the $G$-graded free algebra $F\langle Y, Z \rangle,$
a bilinear form  $\mathfrak{f}_2$, and a linear form $\mathfrak{f}_1$
(see \cite{Svi2}). Here the algebra of graded pure form $*$-polynomials
$\mathcal{S}$ is the free associative commutative algebra
with unit generated by $\mathfrak{f}_2(u_1,u_2),$
$\mathfrak{f}_1(u_3)$  for all nonempty graded $*$-monomials
$u_1, u_2, u_3 \in F\langle Y, Z \rangle.$
Then $FS\langle Y, Z \rangle$ is a $G$-graded algebra
with the grading induced from $F\langle Y, Z \rangle$ assuming
$\deg_{G} f \otimes s = \deg_{G} f,$ for all $f \in
F\langle Y, Z \rangle,$ \ $s \in \mathcal{S}.$
The algebra $FS\langle Y, Z \rangle$ is called free graded
$*$-algebra with forms.

The concept of graded $*$-identities with forms is introduced as usual
with regard to the grading.
Let $A$ be a graded $R$-algebra with involution and forms,
$f(x_{1},\dots,x_{n}) \in FS\langle Y, Z \rangle$
be a graded $*$-polynomial with forms.
$A$ satisfies the graded $*$-identity with forms $f=0$
if $f(a_{1},\dots,a_{n})=0$ holds in $A$ for any $a_{i} \in A$
with $\deg_{\widehat{G}} x_i =\deg_{\widehat{G}} a_i.$
The ideal of graded $*$-identities with forms of an algebra $A$
$\mathrm{SId}^{gi}(A)=\{ f \in FS\langle Y, Z \rangle | A \ \mbox{satisfies} \ f=0 \}$ is a graded
$\mathcal{S}$-ideal of $FS\langle Y, Z \rangle$
invariant with respect to involution and closed
under all graded $*$-endomorphisms \ $\varphi$ of $FS\langle Y, Z \rangle$
which preserve the forms. $\mathrm{SId}^{gi}(A)$ also has the property
that $g_1 \cdot \mathfrak{f}_2(f,g_2),
g_1 \cdot \mathfrak{f}_2(g_2,f), g_1 \cdot \mathfrak{f}_1(f)
\in \mathrm{SId}^{gi}(A)$ \
for any $g_1, g_2 \in FS\langle Y, Z \rangle,$ \
$f \in \mathrm{SId}^{gi}(A).$ Ideals of
$FS\langle Y, Z \rangle$ with all mentioned properties are called $gi$TS-ideals.
Given a $gi$TS-ideal $\widetilde{\Gamma}$
we define the relatively free graded $*$-algebra with forms of infinite rank $\overline{FS}\langle
Y, Z\rangle=FS\langle Y, Z \rangle /\widetilde{\Gamma},$ and of a rank $\nu$ $\overline{FS}\langle
Y_{(\nu)}, Z_{(\nu)}\rangle=FS\langle Y_{(\nu)}, Z_{(\nu)} \rangle /(\widetilde{\Gamma} \cap FS\langle Y_{(\nu)}, Z_{(\nu)} \rangle).$ The equality of graded $*$-polynomials with forms modulo $\widetilde{\Gamma}$ is defined similarly to non-graded case \cite{Svi2}.
We denote
also by $giTS[\mathcal{V}]$ the $gi$TS-ideal generated by a set
$\mathcal{V} \subseteq FS\langle Y, Z \rangle.$

Assume now that $F$ is a field of characteristic zero, and $\sqrt[\mathfrak{m}]{1} \in F.$
Let us define forms on a finite dimensional $G$-graded
$F$-algebra with involution $A=B \oplus J$ with the Jacobson radical $J$,
and the semisimple part $B$.
Consider for any element $b \in B_{\mathfrak{e}}$ the linear operator
$\mathfrak{T}_{b}:B \rightarrow B$  on the graded $*$-subalgebra $B$
defined by
\begin{eqnarray} \label{Oper}
\mathfrak{T}_b(c)= b \circ c, \quad  c \in B.
\end{eqnarray}
It is clear
that $\mathfrak{T}_{\alpha_1 b_1 + \alpha_2 b_2}=
\alpha_1 \mathfrak{T}_{b_1}+\alpha_2 \mathfrak{T}_{b_2}$ for all
$\alpha_i \in F,$ \ $b_i \in B_{\mathfrak{e}}.$ If $b \in
B_{\mathfrak{e}}^{+}$ is
symmetric element then the subspaces $B_{\theta}^{\delta}$
are stable under $\mathfrak{T}_{b}$ for all $(\delta,\theta) \in \widehat{G}.$
If $b \in B_{\mathfrak{e}}^{-}$ is skew-symmetric
then $\mathfrak{T}_{b}(B_{\theta}^{+}) \subseteq B_{\theta}^{-},$
and $\mathfrak{T}_{b}(B_{\theta}^{-}) \subseteq B_{\theta}^{+},$ \
$\theta \in G.$ Particularly, the trace of
the operator $\mathfrak{T}_{b}$ is zero for any
$b \in B_{\mathfrak{e}}^{-}.$

Then the bilinear form $\mathfrak{f}_2:A^2 \rightarrow F,$
and the linear form $\mathfrak{f}_1:A \rightarrow F$
are defined on $A$ by the rules
\begin{eqnarray} \label{Atrace}
&&\mathfrak{f}_2(a_1,a_2)=
\mathfrak{f}_2(b^{(\mathfrak{e})}_1,b^{(\mathfrak{e})}_2)=
\mathrm{Tr}(\mathfrak{T}_{b^{(\mathfrak{e})}_1} \cdot
\mathfrak{T}_{b^{(\mathfrak{e})}_2}),
\nonumber \\
&&\mathfrak{f}_1(a_1)=\mathfrak{f}_1(b^{(\mathfrak{e})}_1)=
\mathrm{Tr}(\mathfrak{T}_{b^{(\mathfrak{e})}_1}), \\
&&a_i=b_i+r_i \in A, \ b_i=\sum_{\theta \in G} b^{(\theta)}_i \in B, \  r \in J, \ \  b^{(\theta)}_i \in B_{\theta}, \ \theta \in G, \nonumber
\end{eqnarray}
where $\mathfrak{T}_{1} \cdot \mathfrak{T}_{2}$ is the product of linear
operators, and $\mathrm{Tr}$ is the usual trace function of linear operator.
Suppose that $A=A_1 \times \dots \times A_\rho.$ Observe that in this case the restrictions on $A_i$ of the
forms $\mathfrak{f}_1,$ $\mathfrak{f}_2$ of $A$ coincide with the forms defined by (\ref{Atrace})
on $A_i$ directly.
It is clear that $\mathfrak{f}_2$ is a symmetric form satisfying
$\mathfrak{f}_2(r,a)=0$ for any $r \in J,$ \  $a \in A,$ \
$\mathfrak{f}_2(a_1,a_2)=0$ for any $a_1 \in A^{-},$ \ $a_2 \in A^{+},$
and $\mathfrak{f}_2(a_1,a_2)=0$ for any $a_1 \in A_{\theta},$
$\theta \ne \mathfrak{e},$\ $a_2 \in A.$
The linear form $\mathfrak{f}_1$ also satisfies
$\mathfrak{f}_1(r)=0$ for any $r \in J,$ \
$\mathfrak{f}_1(a)=0$ for any $a \in A^{-},$ and
$\mathfrak{f}_1(a)=0$ for any $a \in A_{\theta},$
$\theta \ne \mathfrak{e}.$
Particularly, the next lemma holds.

\begin{lemma} \label{Traceid10}
Let $A$ be a finite dimensional graded $F$-algebra with involution and with the forms
(\ref{Atrace}).  Given a graded form $*$-polynomial $h \in FS\langle Y, Z \rangle,$ and
variables $x_1,$ $x_2, x_3 \in Y \cup Z$ with the exception of three cases
$x_1, x_2 \in Y^{\mathfrak{e}},$ or $x_1, x_2 \in Z^{\mathfrak{e}},$ or
\  $x_3 \in Y^{\mathfrak{e}}$ \ $A$ satisfies graded $*$-identities with forms \\
$\begin{array}{llll}
\qquad \qquad
&\mathfrak{f}_1(x_3) \cdot h=0, \qquad
&\mathfrak{f}_2(x_1,x_2) \cdot h=0.
\end{array}$
\end{lemma}

Applying the arguments of Lemma 21 \cite{Svi2} and considering
restrictions of the corresponding operators on $B^{\delta}_{\theta}$
($(\delta,\theta) \in \widehat{G}$) we obtain the following lemma in graded case.
Observe that here it is enough to consider semisimple or radical evaluations of variables
(not necessary elementary ones).

\begin{lemma} \label{Traceid1}
Given a finite dimensional $G$-graded
$*$-algebra $A$ with the forms (\ref{Atrace}) over a field $F,$
and a graded $*$-polynomial $f \in F\langle Y, Z \rangle$
of type $(\mathrm{dims}_{gi} A,\mathrm{nd}(A)-1,1)$
suppose that $\{ x_1, \dots, x_t \} \in Y \cup Z$ are variables on which $f$ is $(\mathrm{dims}_{gi} A)$-alternating ($t=\dim B$).
Then $A$ satisfies the graded $*$-identities with forms \\
$\begin{array}{ll}
\qquad
&\mathfrak{f}_2(y_1,y_2)f=
\sum_{i=1}^{t} f|_{x_i:=y_1 \circ (y_2 \circ x_i)}, \quad
y_1, y_2 \in Y^{\mathfrak{e}}, \\
\qquad
&\mathfrak{f}_2(z_1,z_2)f=
\sum_{i=1}^{t} f|_{x_i:=z_1 \circ (z_2 \circ x_i)}, \quad
z_1, z_2 \in Z^{\mathfrak{e}}, \\
\qquad
&\mathfrak{f}_1(y)f=
\sum_{i=1}^{t} f|_{x_i:=y \circ x_i}, \quad
y \in Y^{\mathfrak{e}}. \\
\end{array}$
\end{lemma}

\begin{lemma} \label{Traceid2} \qquad
Let $f(\widetilde{x}_1,\dots,\widetilde{x}_k) \in F\langle Y, Z \rangle$
be a graded $*$-polynomial of a type
$(\beta;\gamma-1;1)$ (for some $\beta \in \mathbb{N}_0^{2 \mathfrak{m}},$ \ $\gamma \in \mathbb{N}$), and $s(\zeta_1,\dots,\zeta_d) \in \mathcal{S}$ a graded pure form $*$-polynomial ($\{\zeta_1,\dots,\zeta_d\}
\subseteq Y \cup Z$). Then
there exists a graded $*$-polynomial
$g_s(\widetilde{x}_1,\dots,\widetilde{x}_k,\zeta_1,\dots,\zeta_d) \in giT[f]$
such that any finite
dimensional $G$-graded $*$-algebra $A$
with forms (\ref{Atrace}) having parameter
$\mathrm{par}_{gi}(A)=(\beta;\gamma)$ satisfies the graded $*$-identity
with forms
\begin{equation*}
s(\zeta_1,\dots,\zeta_d) \cdot f(\widetilde{x}_1,\dots,\widetilde{x}_k) -
g_s(\widetilde{x}_1,\dots,\widetilde{x}_k,\zeta_1,\dots,\zeta_d) =0.
\end{equation*}
\end{lemma}
\noindent {\bf Proof.} Assume that
$f$ is $(\mathrm{dims}_{gi} A)$-alternating on
$\{ \widetilde{x}_1,\dots,\widetilde{x}_t \},$ \ $t=\dim B.$
Suppose that $w_i$ are $\widehat{G}$-homogeneous
polynomials of $\widehat{G}$-degree  $\deg_{\widehat{G}} \widetilde{x}_i$
($i=1, \dots, t$), and
$\tilde{\zeta}_{j} \in Y^{\mathfrak{e}} \cup Z^{\mathfrak{e}}$ are variables
satisfying the claims of Lemma \ref{Traceid1}.
Applying consequently Lemma \ref{Traceid1} we obtain that
$\mathfrak{f}_2(\tilde{\zeta}_{1},\tilde{\zeta}_{2}) \cdots
\mathfrak{f}_2(\tilde{\zeta}_{2 n_2-1},\tilde{\zeta}_{2 n_2}) \cdot
\mathfrak{f}_1(\tilde{\zeta}_{2 n_2+1}) \cdots
\mathfrak{f}_1(\tilde{\zeta}_{2 n_2+n_1}) \times$
$f(w_1,\dots,w_{t},\widetilde{X})=\sum_{l=1}^{\widetilde{n}}
f(\widetilde{w}_{l 1},\dots,\widetilde{w}_{l t},\widetilde{X})$
$(\mathrm{mod} \ \mathrm{SId}^{gi}(A))$  for some
$\widehat{G}$-homogeneous polynomials $\widetilde{w}_{l i}$ such that
$\deg_{\widehat{G}} \widetilde{w}_{l i} =\deg_{\widehat{G}} \widetilde{x}_i$
for all $l$ ($i=1, \dots, t$). In fact it is sufficient to consider as $w_i$
right normed jordan monomials of the form
$(\tilde{\zeta}_{j_1} \circ (\tilde{\zeta}_{j_2} \circ (\dots (\tilde{\zeta}_{j_r} \circ  \widetilde{x}_i)))).$ Therefore $\widetilde{w}_{l i}$ are right normed jordan monomials of the
same type.
Replacing $\tilde{\zeta}_{j}$ by homogeneous elements of $F\langle Y,Z \rangle$
of the corresponding $\widehat{G}$-degrees, and applying Lemma \ref{Traceid10}
as in the proof of Lemma 22 \cite{Svi2}
we obtain that
$s(\zeta_1,\dots,\zeta_d) \cdot
f(\widetilde{x}_1,\dots,\widetilde{x}_{k}) = g_s(\widetilde{x}_1,\dots,
\widetilde{x}_{k},\zeta_1,\dots,\zeta_d)
(\mathrm{mod} \ \mathrm{SId}^{gi}(A)).$ Observe that the graded $*$-polynomial
$g_s$ does not depend on $A$ and $g_s \in giT[f].$
\hfill $\Box$

Assume that $A$ is a $G$-graded finite
dimensional $*$-algebra with the Jacobson radical $J,$ and the semisimple part $B.$
Let us denote $t_{\vartheta}=\dim B^{\delta}_{\theta},$ \
$q_{\vartheta}=\dim J^{\delta}_{\theta}$ for any $\vartheta=(\delta,\theta) \in \widehat{G},$
and $t=\sum_{\vartheta \in \widehat{G}} t_\vartheta=\dim B.$
Given a positive integer $\nu$ take
$\Lambda_{\nu} = \{ \lambda_{\vartheta i j} | \vartheta \in \widehat{G}; 1 \le i \le
\nu; 1 \le j \le t_\vartheta+q_\vartheta \},$ and the free commutative associative unitary algebra
$F[\Lambda_{\nu}]^{\#}$ generated by the set $\Lambda_{\nu}.$
Let us consider the extension
$\mathcal{P}_{\nu}(A)=F[\Lambda_{\nu}]^{\#} \otimes_F A$ of $A$ by $F[\Lambda_{\nu}]^{\#}.$

$\mathcal{P}_{\nu}(A)$ is a graded algebra with the involution defined
by $(f \otimes a)^*=f \otimes a^*$
($f \in F[\Lambda_{\nu}]^{\#},$ \ $a \in A$), and the grading
$(\mathcal{P}_{\nu}(A))_\theta=F[\Lambda_{\nu}]^{\#} \otimes_F A_{\theta},$
$\theta \in G.$
The forms $\mathfrak{f}_2,$ $\mathfrak{f}_1$  of $A$ defined by (\ref{Atrace})can be naturally
extended to the $F[\Lambda_{\nu}]^{\#}$-bilinear form
$\mathfrak{f}_2: \mathcal{P}_{\nu}(A)^2 \rightarrow
F[\Lambda_{\nu}]^{\#},$ and $F[\Lambda_{\nu}]^{\#}$-linear
form $\mathfrak{f}_1: \mathcal{P}_{\nu}(A) \rightarrow
F[\Lambda_{\nu}]^{\#}$ respectively.

We call by a Cayley-Hamilton type graded $*$-polynomial a degree homogeneous graded $*$-polynomial with forms of the following type
\begin{eqnarray*}
x^{n}+\sum \limits_{\mathop{\ i_0 + i_1 + \dots + j_{k_2} =n,}
\limits_{\scriptstyle 0 < i_0 < n, \  1 \le k_2+k_1}}
\alpha_{(i),(j)} \ \ x^{i_0} \mathfrak{f}_2(x^{i_1},x^{j_1}) \cdots
\mathfrak{f}_2(x^{i_{k_2}},x^{j_{k_2}})
\mathfrak{f}_1(x^{i_{k_2+1}}) \cdots \mathfrak{f}_1(x^{i_{k_2+k_1}}),
\end{eqnarray*}
where $\alpha_{(i),(j)} \in F,$ \ $x=y+z,$ \ $y \in Y^{\mathfrak{e}},$ \
$z \in Z^{\mathfrak{e}}.$  Note that here $i_l, j_l > 0$ ($l \ge 0$).
A Cayley-Hamilton type polynomial is not $\widehat{G}$-homogeneous, but it
is $G$-homogeneous of the neutral degree.

\begin{lemma} \label{HamKel}
$\mathcal{P}_{\nu}(A)$ satisfies a Cayley-Hamilton type graded $*$-identity
$\mathcal{K}_{3t+1}^{\mathrm{nd}(A)}(x)=0$ for some
Cayley-Hamilton type graded $*$-polynomial $\mathcal{K}_{3t+1}(x)$
of degree $3t+1,$ \ $t=\dim B,$ \ $x=y+z,$ \ \ $y \in Y^{\mathfrak{e}},$ \
$z \in Z^{\mathfrak{e}}.$
\end{lemma}
\noindent {\bf Proof.}
By Lemma 23 \cite{Svi2} the algebra $\mathcal{P}_{\nu}(A)$ satisfies the
non-graded $*$-identity $\mathcal{K}_{3t+1}^{\mathrm{nd}(A)}(x)=0,$ where
$\mathcal{K}_{3t+1}(x)$ is a Cayley-Hamilton type non-graded $*$-polynomial
of degree $3t+1$ with the forms defined by (15), (16) in \cite{Svi2}. Particularly, $\mathcal{K}_{3t+1}^{\mathrm{nd}(A)}(x)=0$ holds for any $x \in (\mathcal{P}_{\nu}(A))_{\mathfrak{e}}.$
Observe that
for all powers of an element $x \in (\mathcal{P}_{\nu}(A))_{\mathfrak{e}}$ the definition (15) of the forms
$\mathfrak{f}_1,$ $\mathfrak{f}_2$ in non-graded case given in \cite{Svi2}
coincides with the corresponding definition
(\ref{Atrace}) in the $G$-graded case.
\hfill $\Box$

Let  $\{\hat{b}_{\vartheta 1}, \dots,
\hat{b}_{\vartheta t_\vartheta}\}$ be a basis of the $\widehat{G}$-homogeneous part
$B_{\theta}^{\delta}$ of a semisimple part $B$ of $A,$ and
$\{\hat{r}_{\vartheta 1}, \dots, \hat{r}_{\vartheta
q_{\vartheta}}\}$ a basis of the $\widehat{G}$-homogeneous part $J_{\theta}^{\delta}$ of the
Jacobson radical $J=J(A)$ of $A,$
\ $\vartheta =(\delta,\theta) \in \widehat{G}.$
If $A$ has an elementary decomposition then all these bases can be chosen in the set
$\bigcup_{\vartheta \in \widehat{G}} \bigl( D_{\vartheta} \cup U_{\vartheta} \bigr)$
((\ref{basisD}), (\ref{basisU}), Lemma \ref{Pierce}). Let us take
the elements
\begin{eqnarray} \label{genset}
\mathfrak{y}_{\vartheta i}=\sum_{j=1}^{t_\vartheta} \lambda_{\vartheta i j} \otimes
\hat{b}_{\vartheta j} + \sum_{j=1}^{q_\vartheta} \lambda_{\vartheta i
j+t_\vartheta} \otimes \hat{r}_{\vartheta j} \  \in  \mathcal{P}_{\nu}(A),
\qquad \vartheta \in \widehat{G}, \  \ 1 \le i \le \nu.
\end{eqnarray}
All elements $\mathfrak{y}_{\vartheta i}$ are $\widehat{G}$-homogeneous of
$\widehat{G}$-degree $\vartheta.$ Denote by $\mathfrak{Y}_{\nu}=\{ \mathfrak{y}_{\vartheta i} |
\vartheta \in \widehat{G}; \  1 \le i \le \nu \}$ the set of these elements.
Consider the $G$-graded $*$-invariant
$F$-subalgebra $\mathcal{F}_{\nu}(A)$ of $\mathcal{P}_{\nu}(A)$
generated by $\mathfrak{Y}_{\nu}$.
Consider any map $\varphi$ of the generators $\mathfrak{Y}_{\nu}$ to
arbitrary $\widehat{G}$-homogeneous elements $\widetilde{a}_{\vartheta i} \in A$
of the corresponding $\widehat{G}$-degrees
\begin{equation} \label{hom1}
\varphi:\mathfrak{y}_{\vartheta i} \mapsto \widetilde{a}_{\vartheta
i}=\sum_{j=1}^{t_\vartheta} \widetilde{\alpha}_{\vartheta i j}
\hat{b}_{\vartheta j} + \sum_{j=1}^{q_\vartheta}
\widetilde{\alpha}_{\vartheta i j+t_\vartheta} \hat{r}_{\vartheta j}
\quad (\vartheta \in \widehat{G}; \quad i=1,\dots,\nu),
\end{equation}
here $\widetilde{\alpha}_{\vartheta i j} \in F$.
It is clear that $\varphi$
can be extended to the graded $*$-homomorphism of $F$-algebras
$\varphi:\mathcal{F}_{\nu}(A) \rightarrow A$, also inducing the
graded $*$-homomorphism $\widetilde{\varphi}:\mathcal{P}_{\nu}(A)
\rightarrow A$ by the equalities
\begin{eqnarray} \label{hom2}
\widetilde{\varphi}((\lambda_{\vartheta_1 i_1 j_1} \cdots
\lambda_{\vartheta_k i_k j_k})\otimes a)=(\widetilde{\alpha}_{\vartheta_1
i_1 j_1} \cdots \widetilde{\alpha}_{\vartheta_k i_k j_k}) \cdot a
\qquad \forall a \in A.
\end{eqnarray}
The graded $*$-homomorphism $\widetilde{\varphi}$ preserves the
forms on $\mathcal{P}_{\nu}(A)$ and $A$ defined by (\ref{Atrace}).

Elements of $\mathcal{F}_{\nu}(A)$
are called quasi-polynomials on the variables $\mathfrak{Y}_{\nu}.$ Products
of the generators $\mathfrak{y}_{\vartheta i} \in \mathfrak{Y}_{\nu}$ of the algebra
$\mathcal{F}_{\nu}(A)$ are called quasi-monomials. We have
also $\mathrm{Id}^{gi}(\mathcal{F}_{\nu}(A)) \supseteq
\mathrm{Id}^{gi}(\mathcal{P}_{\nu}(A))=\mathrm{Id}^{gi}(A)$ for any $\nu
\in \mathbb{N}.$

$\mathcal{F}_{\nu}(A)$ is a finitely generated PI-algebra. By the
Shirshov's theorem on height \cite{Shirsh} $\mathcal{F}_{\nu}(A)$
has a finite height and a finite Shirshov's basis that can be chosen
in the set of monomials over the generators (\cite{BelRow}, \cite{Shirsh}).
More precisely there exist a natural
$\mathcal{H},$ and quasi-monomials $w_1,\dots,w_d \in \mathcal{F}_{\nu}(A)$ such that
any element $u \in \mathcal{F}_{\nu}(A)$ has the form
$u=\sum_{(i)=(i_1,\dots,i_k)}
 \alpha_{(i)} \  w_{i_1}^{c_1} \dots w_{i_k}^{c_k},$
where $k \le \mathcal{H},$ \  $\{ i_1,\dots,i_k \} \subseteq
\{1,\dots,d\},$ \  $c_j \in \mathbb{N},$ \  $\alpha_{(i)} \in F.$
Observe that the elements $w_i$ are $G$-homogeneous, but not necessarily
$\widehat{G}$-homogeneous.

Consider the polynomials
$\hat{\mathfrak{s}}_{i,(l_1,l_2)}=
\mathfrak{f}_2(w_i^{\mathfrak{m} l_1},w_i^{\mathfrak{m} l_2}) \in
F[\Lambda_{\nu}]^{\#}$ ($i=1,\dots,d,$ \
$l_1,l_2=1,\dots,3t,$ \ $t=\dim B$),
and
$\hat{\mathfrak{s}}_{i,l}=
\mathfrak{f}_1(w_i^{\mathfrak{m} l}) \in
F[\Lambda_{\nu}]^{\#}$
($i=1,\dots,d,$ \ $l=1,\dots,3t$).
Then
$\widehat{F}=F[
\hat{\mathfrak{s}}_{i,(l_1,l_2)},
\hat{\mathfrak{s}}_{i, l_1} \ | \ 1 \le i \le d; \ 1
\le l_1, l_2 \le 3t \ ]^{\#}$ is
the associative commutative unitary $F$-subalgebra
of $F[\Lambda_{\nu}]^{\#}$
generated by $\{ \hat{\mathfrak{s}}_{i,(l_1,l_2)},
\hat{\mathfrak{s}}_{i, l} \}$, and
the unit of $F[\Lambda_{\nu}]^{\#}.$

Take the graded $*$-invariant $\widehat{F}$-subalgebra
$\mathcal{T}_{\nu}(A)=\widehat{F} \mathcal{F}_{\nu}(A)$ of
$\mathcal{P}_{\nu}(A).$ Then $\mathcal{F}_{\nu}(A)$ is
a graded $*$-subalgebra of $\mathcal{T}_{\nu}(A).$ Arbitrary
map of the form (\ref{hom1}) can be properly and uniquely extended to
the graded $*$-homomorphism from $\mathcal{T}_{\nu}(A)$ to $A$ preserving the
forms (it is the restriction on $\mathcal{T}_{\nu}(A)$ of
$\widetilde{\varphi}$ defined by (\ref{hom2})).
Since $w_i^{\mathfrak{m}} \in (\mathcal{P}_{\nu}(A))_{\mathfrak{e}}$ \ ($i=1,\dots,d$)
then by Lemma \ref{HamKel} all elements $w_i^{\mathfrak{m}}$
are algebraic of the degree $\mathrm{nd}(A)(3t+1)$ over
$\widehat{F}.$
Therefore by Shirshov's theorem on height
$\mathcal{T}_{\nu}(A)$ is finitely generated
$\widehat{F}$-module, where $\widehat{F}$ is Noetherian. By
theorem of Beidar \cite{Beid} the algebra $\mathcal{T}_{\nu}(A)$
is representable.

Observe that all elements of the set $\mathfrak{Y}_{\nu}$
are $\widehat{G}$-homogeneous and  homomorphisms defined by (\ref{hom1}), (\ref{hom2})
are graded. Therefore we obtain the graded versions of Remark 2, and Lemmas 24, 25 \cite{Svi2}
repeating their arguments and using
the graded versions of Lemma 22, and Lemma 5
of \cite{Svi2} (Lemma \ref{Traceid2}, and Lemma \ref{Repr} respectively).
As consequences we have also the $G$-graded versions of Lemmas 26, 27 \cite{Svi2}.
In the graded versions of Lemmas 25-27 we assume that Assumption \ref{Class} is true over any
algebraically closed extension of $F.$ Besides that the graded versions of Lemmas 26, 27
are proved for the $gi$T-ideal $\Gamma$ of a finitely generated $G$-graded PI-algebra with
involution, i. e. we assume that $\Gamma$ contains a non-trivial T-ideal.

\section{Graded $*$-identities of finitely generated algebras.}

Now we can state the relation between graded $*$-identities
of finitely generated PI-algebras and graded $*$-identities
of finite dimensional algebras
assuming that Assumption \ref{Class} is true in our case.

\begin{theorem} \label{*PI-gen}
Let $G$ be a finite abelian group, and $F$ a field of characteristic zero
containing a primitive root of $1$ of the degree
$\mathfrak{m}=|G|.$ Suppose that Assumption \ref{Class} holds for the group $G$
over any algebraically closed extension $\widetilde{F}$ of $F.$
Consider a finitely generated $G$-graded associative PI-algebra $D$ over $F$
with graded involution. Then
the $gi$T-ideal of graded $*$-identities of $D$ coincides
with the $gi$T-ideal of graded $*$-identities of some finite dimensional
$G$-graded associative $F$-algebra with graded involution.
\end{theorem}
\noindent {\bf Proof.} Let $\Gamma = Id^{gi}(D).$  We use the induction on the
Kemer index $\mathrm{ind}_{gi}(\Gamma)=\kappa=(\beta;\gamma)$ of $\Gamma.$

\underline{The base of induction.} Let
$\mathrm{ind}_{gi}(\Gamma)=(\beta;\gamma)$ with $\beta=(0,\dots,0).$
Then $D$ is a nilpotent finitely generated algebra. Hence, $D$ is finite
dimensional.

\underline{The inductive step.} The $G$-graded versions of
Lemmas 11, 26, 27 \cite{Svi2} imply that
$\Gamma \supseteq \mathrm{Id}^{gi}(A),$ where
$A=\mathcal{O}(A) \times \mathcal{Y}(A)$ is a finite dimensional $G$-graded algebra with
involution, such that
$\mathrm{ind}_{gi}(\Gamma)=\mathrm{ind}_{gi}(A)=\kappa.$ Moreover,
$S_{\tilde{\mu}}(\Gamma)=S_{\tilde{\mu}}(\mathcal{O}(A))=S_{\tilde{\mu}}(A)
\subseteq \mathrm{Id}^{gi}(\mathcal{Y}(A))$ for some $\tilde{\mu} \in
\mathbb{N}.$ Here $\mathcal{O}(A),$ $\mathcal{Y}(A)$ are finite dimensional $G$-graded $*$-algebras
with elementary decomposition. $\mathcal{O}(A)=A_1 \times \cdots \times A_{\rho},$
where $A_i$ are $gi$-reduced algebras, $\mathrm{ind}_{gi}(\mathcal{O}(A))=\mathrm{ind}_{gi}(A_i)=\kappa$
\ $\forall i=1,\dots,\rho,$ and $\mathrm{ind}_{gi}(\mathcal{Y}(A)) < \kappa$ (see the graded version
of Lemma 11 \cite{Svi2}).

Denote  $(t_1,\dots,t_{2 \mathfrak{m}})=\beta(\Gamma)=\mathrm{dims}_{gi} A_i,$
\ $t=\sum_{j=1}^{2 \mathfrak{m}} t_j;$ \ $\gamma=\gamma(\Gamma)=\mathrm{nd}(A_i)$
(for all $i=1,\dots,\rho$). Let us take for any $i=1,\dots,\rho$
the algebra $\widetilde{A}_i=\mathcal{R}_{q_i,s}(A_i)$ defined by
(\ref{FRad}) for $A_i$ with $q_i=\dim_F A_i,$
\  $s=(t+1)(\gamma+\tilde{\mu}).$ $\widetilde{A}_i$ is
a finite dimensional graded $*$-algebra with elementary decomposition.
We have also $\Gamma_i=\mathrm{Id}^{gi}(\widetilde{A}_i)=$ $\mathrm{Id}^{gi}(A_i),$
and $\mathrm{dims}_{gi} \widetilde{A}_i=$ $\mathrm{dims}_{gi} A_i=$ $\beta.$
The Jacobson radical $J(\widetilde{A}_i)=(Y_{(q_i)},Z_{(q_i)})/I$ is
nilpotent of the degree at most $s=(t+1)(\gamma+\tilde{\mu}),$
where $I=\Gamma_i(B_i(Y_{(q_i)},Z_{(q_i)}))+(Y_{(q_i)},Z_{(q_i)})^s.$ Here the algebra
$B_i$ can be considered as the semisimple part of $A_i$ and of
$\widetilde{A}_i$ simultaneously (the graded version of Lemma 3 \cite{Svi2}).
Take $\widetilde{A}=\times_{i=1}^\rho
\widetilde{A}_i,$ and $\nu=\mathrm{rkh}(D).$
By the graded version of Remark 2 \cite{Svi2}, and Lemma \ref{Repr} there
exists an $F$-finite dimensional $G$-graded algebra with graded
involution and elementary decomposition $C$ such that
$\mathrm{Id}^{gi}(C)=\mathrm{Id}^{gi}(\mathcal{T}_{\nu}(\widetilde{A})/\Gamma
(\mathcal{T}_{\nu}(\widetilde{A}))).$

Let us denote $\widetilde{D}_{\nu}=F\langle Y_{(\nu)}, Z_{(\nu)}
\rangle/\bigl((\Gamma + K_{\widetilde{\mu}}(\Gamma)) \cap
F\langle Y_{(\nu)}, Z_{(\nu)} \rangle \bigr).$ The graded versions of Lemmas 6, 8 \cite{Svi2}
imply that
$\mathrm{ind}_{gi}(\widetilde{D}_{\nu}) \le \mathrm{ind}_{gi}(\Gamma +
K_{\tilde{\mu}}(\Gamma)) < \mathrm{ind}_{gi}(\Gamma)$. By the inductive
step we obtain
$\mathrm{Id}^{gi}(\widetilde{D}_{\nu})=\mathrm{Id}^{gi}(\widetilde{U}),$
where $\widetilde{U}$ is a finite dimensional over $F$ $G$-graded $*$-algebra.
The graded version of Remark 2 \cite{Svi2} yields $\Gamma
\subseteq \mathrm{Id}^{gi}(C \times \widetilde{U}).$

Consider a multilinear polynomial $f(\tilde{x}_1,\dots,\tilde{x}_m) \in
\mathrm{Id}^{gi}(C \times \widetilde{U})$ in variables
$\tilde{x}_j \in Y \cup Z.$
Let us take any multihomogeneous with respect to degrees of variables
and $\widehat{G}$-homogeneous $*$-polynomials
$w_1, \dots, w_m \in F\langle Y_{(\nu)}, Z_{(\nu)} \rangle$
($\deg_{\widehat{G}} w_j = \deg_{\widehat{G}} \tilde{x}_j,$ \ $j=1,\dots,m$).
We have $f(w_1,\dots,w_m)=g+h$ for some multihomogeneous graded
$*$-polynomials $g \in \Gamma,$ \ $h \in K_{\tilde{\mu}}(\Gamma)$
also depending on $Y_{(\nu)} \cup Z_{(\nu)}.$
Then by the graded version of Lemma 24 \cite{Svi2} we obtain $h=f(w_1,\dots,w_m)-g \in \mathcal{S}
\Gamma + \mathrm{SId}^{gi}(\widetilde{A}_i)$ for any
$i=1,\dots,\rho.$ Hence
$\tilde{h}(\tilde{x}_1,\dots,\tilde{x}_n) \in \mathcal{S} \Gamma +
\mathrm{SId}^{gi}(\widetilde{A}_i)$ holds also for the full
linearization $\tilde{h}$ of $h$.
The graded version of Lemma 18 \cite{Svi2} implies that $\tilde{h}$ is
exact for $A_i$ ($i=1,\dots,\rho$).

Fix any $i=1,\dots,\rho.$ Assume that $(c_1,\dots,c_{q_i})$ is a basis of $A_i$ consisting of $\widehat{G}$-homogeneous elements
chosen in $D \cup U$ (Lemma \ref{Pierce}), and fix the order of the basic elements.
Suppose that $\bar{a}$ is
an elementary complete evaluation of $\tilde{h}$ in the algebra
$A_i$ with $\gamma-1$ radical elements.
By Lemma \ref{Gammasub} there exists a polynomial
$h_{\tilde{\mu}}(\mathcal{Z}_1,\dots,\mathcal{Z}_{\gamma-1+\tilde{\mu}},
\mathcal{X}) \in \mathcal{S} \Gamma +
\mathrm{SId}^{gi}(\widetilde{A}_i)$ of type $(\beta,\gamma-1,\tilde{\mu})$, and
an elementary evaluation $\bar{u}$ of $h_{\tilde{\mu}}$ in $A_i$
such that $h_{\tilde{\mu}}(\bar{u})=\tilde{h}(\bar{a}).$

Moreover, $h_{\tilde{\mu}}$ is alternating in any
$\mathcal{Z}_j$ ($j=1,\dots,\gamma-1+\tilde{\mu}$), and
all variables from $\mathcal{X}$
are replaced by semisimple elements.
Then we have
\begin{eqnarray} \label{lhdec}
\alpha_2 h_{\tilde{\mu}}=\bigl(
\prod_{m=1}^{\gamma-1+\tilde{\mu}}
\prod_{\vartheta \in \widehat{G}}
\mathcal{A}_{\mathcal{Z}^{\vartheta}_{m}}
 \bigr) h_{\tilde{\mu}}=
\sum_j  \bigl( \prod_{m=1}^{\gamma-1+\tilde{\mu}}
\prod_{\vartheta \in \widehat{G}}
\mathcal{A}_{\mathcal{Z}^{\vartheta}_{m}}
\bigr) \bigl(
\mathfrak{\tilde{s}}_j \tilde{g}_j \bigr) (\mathrm{mod } \
\mathrm{SId}^{gi}(\widetilde{A}_i)),
\end{eqnarray}
where $\mathcal{Z}^{\vartheta}_m$ is the subset of $\widehat{G}$-homogeneous variables of
$\mathcal{Z}_m$ of complete degree $\vartheta=(\delta,\theta) \in \widehat{G};$
\  $\alpha_2 \in F,$ $\alpha_2 \ne 0;$ \ $\tilde{g}_j
\in \Gamma,$ \ $\mathfrak{\tilde{s}}_j \in \mathcal{S}.$
Denote
by $\{\zeta_1,\dots,\zeta_{\hat{n}} \}$
the variables $\mathcal{Z} \cup
\mathcal{X}$ of $h_{\tilde{\mu}}$ (the first $(t+1)(\gamma-1)+t \tilde{\mu}$ variables are
from $\mathcal{Z} = \bigcup_{m=1}^{\gamma-1+\tilde{\mu}}
\mathcal{Z}_m$), and by $\mathcal{Z}^{\vartheta} =
\bigcup_{m=1}^{\gamma-1+\tilde{\mu}} \mathcal{Z}^{\vartheta}_m$
the $\widehat{G}$-homogeneous part of variables $\mathcal{Z}$ of complete
degree $\vartheta \in \widehat{G}.$

Let $u_k \in A_i$ be an element of the mentioned evaluation $\bar{u}=(u_1,\dots,u_{\hat{n}})$ of $h_{\tilde{\mu}}.$
We take in the algebra $\widetilde{A}_i$ the elements
$\bar{y}_{\pi(k) \theta}=y_{\pi(k) \theta}+I,$ \ $\bar{z}_{\pi(k) \theta}=z_{\pi(k) \theta}+I,$
where $y_{\pi(k) \theta} \in Y_{(q_i)}^{\theta},$ \ $z_{\pi(k) \theta} \in Z_{(q_i)}^{\theta}$
are variables, and $u_k=c_{\pi(k)}$ ($\pi(k)$ is the ordinal number of
the element $u_k$ in our basis of $A_i,$ $1 \le \pi(k) \le q_i$), $\deg_{G} u_k =\theta.$

Consider in the algebra $\widetilde{A}_i$ the following evaluation of
$h_{\tilde{\mu}}(\zeta_1,\dots,\zeta_{\hat{n}})$
\begin{eqnarray} \label{Ai}
&& \zeta_k=\bar{y}_{\pi(k) \theta} \in J(\widetilde{A}_i) \qquad
\mbox{ if } \
\zeta_k \in \mathcal{Z}^{(+,\theta)} \  (\theta \in G), \nonumber \\
&& \zeta_k=\bar{z}_{\pi(k) \theta} \in J(\widetilde{A}_i) \qquad
\mbox{ if } \
\zeta_k \in \mathcal{Z}^{(-,\theta)} \  (\theta \in G), \nonumber \\
&& \zeta_k=u_k
\qquad \mbox{ if } \ \zeta_k \in \mathcal{X}.
\end{eqnarray}

If a graded pure form polynomial $\mathfrak{\tilde{s}}_j$ in (\ref{lhdec})
depends essentially on $\mathcal{Z}$ then $\bigl(
\prod_{m=1}^{\gamma-1+\tilde{\mu}}
\prod_{\vartheta \in \widehat{G}}
\mathcal{A}_{\mathcal{Z}^{\vartheta}_{m}}
\bigr) \bigl(
\mathfrak{\tilde{s}}_j \tilde{g}_j \bigr)|_{(\ref{Ai})}=0,$
since the forms are zero on radical
elements (see (\ref{Atrace})). If $\mathfrak{\tilde{s}}_j$ does not depend
on $\mathcal{Z}$ then $\bigl(
\prod_{m=1}^{\gamma-1+\tilde{\mu}}
\prod_{\vartheta \in \widehat{G}}
\mathcal{A}_{\mathcal{Z}^{\vartheta}_{m}}
\bigr) \bigl(
\mathfrak{\tilde{s}}_j \tilde{g}_j \bigr) =$
$\mathfrak{\tilde{s}}_j \tilde{\tilde{g}}_j,$ where
$\tilde{\tilde{g}}_j=\bigl(
\prod_{m=1}^{\gamma-1+\tilde{\mu}}
\prod_{\vartheta \in \widehat{G}}
\mathcal{A}_{\mathcal{Z}^{\vartheta}_{m}}
\bigr)
\tilde{g}_j \in \Gamma.$
If $\tilde{\tilde{g}}_j|_{(\ref{Ai})} \ne 0$ in $\widetilde{A}_i$ then one
of the multihomogeneous on degrees components of
$\tilde{\tilde{g}}_j$ is a $\tilde{\mu}$-boundary polynomial for
$\widetilde{A}_i.$ And it is not a $\tilde{\mu}$-boundary polynomial for $\Gamma,$
because it belongs to $\Gamma.$ It implies $S_{\tilde{\mu}}(A) \ne S_{\tilde{\mu}}(\Gamma),$ that contradicts to the properties of $A.$ Therefore
$\tilde{\tilde{g}}_j|_{(\ref{Ai})}=0.$ Thus in any case
$h_{\tilde{\mu}}|_{(\ref{Ai})}=0$ holds in the algebra
$\widetilde{A}_i.$ Consider in the algebra $B_i(Y_{(q_i)},Z_{(q_i)})$ the elements
\begin{eqnarray*}
&&v_k=y_{\pi(k) \theta} \qquad
\mbox{ if } \
\zeta_k \in \mathcal{Z}^{(+,\theta)} \  (\theta \in G),  \\
&& v_k=z_{\pi(k) \theta} \qquad
\mbox{ if } \
\zeta_k \in \mathcal{Z}^{(-,\theta)} \  (\theta \in G),  \\
&& v_k=u_k
\qquad \mbox{ if } \ \zeta_k \in \mathcal{X}.
\end{eqnarray*}
Hence the evaluation $\zeta_k=v_k$ ($k=1,\dots,\hat{n}$)
of the polynomial $h_{\tilde{\mu}}$ is equal to $h_{\tilde{\mu}}(v_1,\dots,v_{\hat{n}}) \in
I=\Gamma_i(B_i(Y_{(q_i)},Z_{(q_i)}))+(Y_{(q_i)},Z_{(q_i)})^s$ in the algebra
$B_i(Y_{(q_i)},Z_{(q_i)}).$ Since $|\mathcal{Z}| < s,$ the polynomial $h_{\tilde{\mu}}$ is linear
on variables $\mathcal{Z},$ and variables of $\mathcal{X}$ are replaced by semisimple elements
then we obtain $h_{\tilde{\mu}}(v_1,\dots,v_{\hat{n}}) \in
\Gamma_i(B_i(Y_{(q_i)},Z_{(q_i)})).$

Consider the map $\varphi:y_{j \theta} \mapsto
c_{j}$ if $\deg_{\widehat{G}} c_j =(+,\theta),$
and $\varphi:z_{j \theta} \mapsto
c_{j}$ if $\deg_{\widehat{G}} c_j =(-,\theta),$\ $j=1,\dots,q_i$.
It is clear that $\varphi$ can be extended to a
graded $*$-homomorphism $\varphi:B_i(Y_{(q_i)},Z_{(q_i)}) \rightarrow A_i$
assuming $\varphi(b)=b$ for any $b \in B_i.$ Then
$\varphi(h_{\tilde{\mu}}(v_1,\dots,v_{\hat{n}}))=$
$h_{\tilde{\mu}}(\varphi(v_1),\dots,\varphi(v_{\hat{n}}))=$
$h_{\tilde{\mu}}(\bar{u}) \in$ $\varphi(\Gamma_i(B_i(Y_{(q_i)},Z_{(q_i)}))) \subseteq \Gamma_i(A_i)=(0).$

Therefore $\tilde{h}(\bar{a})=h_{\tilde{\mu}}(\bar{u})=0$ holds in $A_i$ for any elementary
complete evaluation $\bar{a} \in A_i^n$ containing $\gamma-1$
radical elements. Since $\tilde{h}$ is a multilinear
exact polynomial for $A_i,$ and $\gamma=\mathrm{nd}(A_i)$
then $\tilde{h} \in
\mathrm{Id}^{gi}(A_i).$ Hence $h \in \cap_{i=1}^\rho
\mathrm{Id}^{gi}(A_i),$ and $h \in
\mathrm{Id}^{gi}(\mathcal{O}(A) \times
\mathcal{Y}(A))=\mathrm{Id}^{gi}(A) \subseteq \Gamma.$ Thus
we have $f(w_1,\dots,w_m)=g+h \in \Gamma$ for all multihomogeneous
$\widehat{G}$-homogeneous
graded $*$-polynomials $w_1, \dots, w_m \in F\langle Y_{(\nu)}, Z_{(\nu)} \rangle$
of corresponding $\widehat{G}$-degrees.
By the graded version of Remark 1 \cite{Svi2} it
implies $\mathrm{Id}^{gi}(C \times \widetilde{U}) \subseteq \Gamma.$

Therefore $\Gamma=\mathrm{Id}^{gi}(C \times \widetilde{U}).$
Theorem is proved. \hfill $\Box$

Theorem \ref{*PI-gen} can be extended for any field of characteristic zero.

\begin{theorem} \label{*PI-gen2}
Let $G$ be a finite abelian group, and $F$ a field of characteristic zero.
Suppose that Assumption \ref{Class} holds for the group $G$
over any algebraically closed extension $\widetilde{F}$ of $F.$
Let $D$ be a finitely generated $G$-graded associative PI-algebra
over $F$ with graded involution.
Then the $gi$T-ideal of graded $*$-identities of $D$ coincides
with the $gi$T-ideal of graded $*$-identities of some finite dimensional over $F$
$G$-graded associative algebra with graded involution.
\end{theorem}
\noindent {\bf Proof.}
Assume that $F$ does not contain $\mathfrak{j}=\sqrt[\mathfrak{m}]{1}.$
Consider the extension $K=F[\mathfrak{j}]$ of $F$ by $\mathfrak{j}.$
It is clear that any algebraically closed extension $\widetilde{F}$ of $F$
contains also $K.$ Since Assumption \ref{Class} holds for the group $G$
over any algebraically closed extension $\widetilde{F}$ of $F$ then
it is true also over any algebraically closed extension $\widetilde{K}$ of $K.$
Consider the algebra $\bar{D}=D \otimes_F K.$
$\bar{D}$ is a finitely generated $K$-algebra with the $G$-grading
$\bar{D}_{\theta}=D_{\theta} \otimes_F K,$ $\theta \in G.$ The graded
involution on $\bar{D}$ is naturally induced from $D$ by equalities
$(a \otimes \alpha)^*=a^* \otimes \alpha,$ for any
$a \in D,$ \ $\alpha \in K.$ It is clear that $D$ can be considered
as an $F$-subalgebra of $\bar{D},$ and the graded $F$-identities with involution of
$D$ and $\bar{D}$ coincide $Id_F^{gi}(D)=Id_F^{gi}(\bar{D})$.
Particularly, $\bar{D}$ is a PI-algebra ($Id_F^{gi}(D) \subseteq Id_{K}^{gi}(\bar{D})$).
By Theorem \ref{*PI-gen} we obtain that $Id_{K}^{gi}(\bar{D})=Id_{K}^{gi}(C)$ over the field $K$
for some $G$-graded algebra $C$ with graded involution, finite dimensional over $K.$
$C$ can be considered also as an $F$-algebra. And as an $F$-algebra $C$ preserves
the same $G$-grading and involution.
Since $K=F[\mathfrak{j}]$ is the finite
extension of $F$ then $C$ is also finite dimensional over $F.$
It is clear that $Id_{F}^{gi}(C)=Id_{K}^{gi}(C) \cap F\langle Y, Z\rangle.$
Therefore we have
$Id_{F}^{gi}(C)=Id_{F}^{gi}(\bar{D})=Id_F^{gi}(D).$
And the $F$-algebra $C$ is the required finite dimensional algebra.
\hfill $\Box$

Observe that the final result is obtained for any base field of characteristic zero.
The unique restrictions that we have are Assumption \ref{Class} and the requirement
for a $gi$T-ideal to contain a non-trivial T-ideal. The second condition is necessary.
An ideal of group-graded identities of a finitely generated algebra can not contain
a non-trivial ordinary non-graded identity (see, e.g., the comment after Theorem 1 \cite{Svi1}).
And a finite dimensional algebra is always a PI-algebra. We have discussed in Section 1 the
conditions which provide this property for a $gi$T-ideal.

\section{PI-representability of $(Z/qZ)$-graded algebras.}

Suppose that $G$ is a cyclic group of order $q,$
where $q$ is a prime number or $q=4.$ We use the additive notation for the
group $G$ in this case.

Consider the function $\chi: \mathbb{Z}/4 \mathbb{Z} \rightarrow \{0, 1 \}$
defined on the group $\mathbb{Z}/4 \mathbb{Z}$ by the rules $\chi(\bar{0})=\chi(\bar{1})=0,$
$\chi(\bar{2})=\chi(\bar{3})=1.$ The next properties of $\chi$ can be checked directly.

\begin{lemma} \label{Propchi}
\quad $\chi(x)+\chi(y)=\chi(x+y) + 1 \ {\rm mod } \ 2$ \quad if $x, y \in \{ \bar{1}, \bar{3}\},$ \\ and \qquad \qquad \quad
$\chi(x)+\chi(y)=\chi(x+y) \ {\rm mod } \  2$ \quad if $x$ or $y$ is even.
\end{lemma}

Recall that an elementary grading
on the matrix algebra $M_k(\widetilde{F})$ is the $G$-grading
defined by a $k$-tuple $(\theta_1,\dots,\theta_k) \in G^k,$ so that $\deg_G(E_{ij})=
-\theta_i + \theta_j$ for any matrix unit $E_{ij}$ (see, e.g., \cite{BahtZaicSeg},
\cite{BahtZaicSeg2}, \cite{BahtShestZ}, \cite{BahtZaic3}).

We obtain the description of
$*$-graded simple finite dimensional algebras over an algebraically closed field $\widetilde{F}$
for the group $G.$
It is based on the classification of simple $G$-graded algebras given
in Lemma \ref{simple} (Theorem 3 \cite{BahtZaicSeg}).

\begin{theorem}  \label{simGr-pClas}
Let $q$ be a prime number or $q=4,$ and
$G$ a cyclic group of order $q.$ Suppose that
$\widetilde{F}$ is an algebraically closed field of characteristic zero, and $C$
is a $G$-graded finite dimensional $\widetilde{F}$-algebra with graded involution.
Then $C$ is $*$-graded simple if and only if $C$ is isomorphic as a
graded $*$-algebra to one of the algebras of the list:
\begin{enumerate}
\item the direct product $\mathcal{B} \times \mathcal{B}^{op}$ of
a graded simple algebra $\mathcal{B}=M_k(\widetilde{F}[H]),$
and its opposite algebra $\mathcal{B}^{op}$ with the exchange involution $\bar{*},$
where $\widetilde{F}[H]$ is the group algebra of the group $H,$ and
$H$ is the trivial group, $G,$ or $H=\{ \bar{0}, \bar{2} \} \leq \mathbb{Z}/4 \mathbb{Z};$
\item the full matrix algebra $M_k(\widetilde{F})$ with an elementary grading and
an elementary involution;
\item the full matrix algebra $M_k(\widetilde{F}[H])$ over the group algebra $\widetilde{F}[H]$
with the grading induced by the natural grading of $\widetilde{F}[H]$ \
($\deg_G \mathcal{X}_{\theta} \  \eta_{\theta} = \theta$), and involution $(\sum_{\theta \in H} \mathcal{X}_{\theta} \  \eta_{\theta})^*=
\sum_{\theta \in H} \mathcal{X}_{\theta}^{\mathfrak{t}} \  \eta_{\theta},$ where $\mathfrak{t}$ is the transpose or symplectic involution on the matrix algebra $M_k(\widetilde{F}),$ \ $\mathcal{X}_{\theta} \in M_k(\widetilde{F}),$ \ $\theta \in H,$ \ $H$ is a cyclic group;
\item the full matrix algebra $M_k(\widetilde{F}[H])$ over the group algebra $\widetilde{F}[H]$
with the grading induced by the natural grading of $\widetilde{F}[H]$
and involution $(\sum_{\theta \in H} \mathcal{X}_{\theta} \  \eta_{\theta})^*=
\sum_{\theta \in H} (-1)^{\theta} \mathcal{X}_{\theta}^{\mathfrak{t}} \  \eta_{\theta},$ where $\mathfrak{t}$ is the transpose or symplectic involution on the matrix algebra $M_k(\widetilde{F}),$ and $H \cong \mathbb{Z}/2 \mathbb{Z},$ or $H \cong \mathbb{Z}/4 \mathbb{Z};$
\item the full matrix algebra $M_k(\widetilde{F}[H])$ over the group algebra of $H=\{ \bar{0}, \bar{2} \}$ with $(\mathbb{Z}/4 \mathbb{Z})$-grading defined as in Lemma \ref{simple} by a $k$-tuple
$(\theta_1,\dots,\theta_k) \in \{ \bar{0}, \bar{1} \}^k$ and an elementary involution.
\end{enumerate}
\end{theorem}
\noindent {\bf Proof.}
It is clear that all alternatives are $*$-graded simple algebras.
Suppose that $C$ is a $*$-graded simple finite dimensional $\widetilde{F}$-algebra. Then $C$ is a $G$-graded semisimple algebra (Lemma \ref{Pierce0}), and it
contains a $G$-graded simple ideal $\mathcal{B}.$ $\mathcal{B}$ is isomorphic as a $G$-graded algebra to $M_k(\widetilde{F}^{\zeta}[H])$ by Lemma \ref{simple}, where $H$ is a subgroup
of $G,$ and $\zeta:H \times H \rightarrow \widetilde{F}^{*}$ is a 2-cocycle on $H.$
The canonical grading of $M_k(\widetilde{F}^{\zeta}[H])$ is defined by a $k$-tuple
$(\theta_1,\dots,\theta_k) \in G^k,$ so that $\deg_G(E_{ij}
\eta_{\xi})=-\theta_i+\xi+\theta_j.$
It is well-known (see, e.g., \cite{Brown}) that the second cohomologies of a cyclic group
in this case are trivial. Thus $\widetilde{F}^{\zeta}[H]$ is isomorphic as a graded algebra to the group algebra $\widetilde{F}[H]$ of $H,$ where $H$ is one of the group of the list: $\{ \mathfrak{e} \},$ $G,$ $\{ \bar{0}, \bar{2} \} \leq \mathbb{Z}/4 \mathbb{Z}.$
Then either $C=\mathcal{B}$ or $C=\mathcal{B} \times \mathcal{B}^*.$ In the last case $C$ is isomorphic to $\mathcal{B} \times \mathcal{B}^{op}$ with the exchange involution.
The isomorphism is given by $\varphi: a+b \mapsto (a,b^*),$ where $a \in B,$
\ $b \in B^*.$ Hence we obtain the first alternative.

Suppose that $C=\mathcal{B}$ is a $G$-graded simple algebra with graded involution. Thus $C \cong M_k(\widetilde{F}[H]),$ where $H \in \{ \{ \mathfrak{e} \}, \{ \bar{0}, \bar{2} \}, G \}.$
If $H=\{ \mathfrak{e} \}$ then $C \cong M_k(\widetilde{F})$ is the full matrix algebra with an elementary
grading and graded involution. The results of Y.A.Bakhturin, I.P.Shestakov, and  M.V.Zaicev
(\cite{BahtShestZ}, \cite{BahtZaic3}) yields in this case that $C$ is isomorphic as a $*$-graded algebra to $M_k(\widetilde{F})$ with an elementary grading and an elementary involution.

Suppose that $H=G,$ and
$C=M_k(\widetilde{F}[G])$ with the grading defined in Lemma \ref{simple} and graded involution.
Consider the algebra $C'=M_k(\widetilde{F}[G])=\mathrm{Span}_F\{ E_{ij} \widetilde{\eta}_{\xi} | i,j=1,\dots,k, \xi \in G\}$ with the $G$-grading induced by the natural grading of $\widetilde{F}[G].$
The $\widetilde{F}$-linear map $\varphi: E_{ij} \eta_{\xi} \mapsto E_{ij}
\widetilde{\eta}_{\xi -\theta_i +\theta_j}$ is a $G$-graded isomorphism of the algebras $C$ and $C'.$
The involution in $C'$ is induced from $C$ by $\varphi.$

A $G$-homogeneous element of $C'$ of degree $\theta \in G$ has the form
$\mathcal{X}_{\theta} \  \widetilde{\eta}_{\theta}=(\mathcal{X}_{\theta} \   \widetilde{\eta}_{\mathfrak{e}}) \cdot (I \  \widetilde{\eta}_{\theta}),$ where $\mathcal{X}_{\theta} \in M_k(\widetilde{F})$ is a matrix,
$I$ is the identity matrix of order $k,$ \ $\mathfrak{e}$ is the unit of the group $G.$
Observe that the element $\mathcal{X}_{\theta} \   \widetilde{\eta}_\mathfrak{e}$ belongs to the neutral component $C'_{\mathfrak{e}}$ of $C'.$  $C'_{\mathfrak{e}} \cong M_k(\widetilde{F}),$ and it is a $*$-invariant subalgebra of $C'.$ By Theorem 4.6.12 \cite{BMM} (see also the proof of Theorem 3.6.8 \cite{GZbook})
the restriction of the involution on $C'_{\mathfrak{e}}$ can be taken as the transpose or
symplectic involution up to an inner automorphism of $C'_{\mathfrak{e}}.$

Consider a generator $\xi$ of the group $G.$ Observe that $I \widetilde{\eta}_{\xi}$ is a central element of $C'$ of $G$-degree $\xi.$ Then $(I \widetilde{\eta}_{\xi})^*$ has the same $G$-degree, and also belongs to the center of $C'.$ Thus $(I \widetilde{\eta}_{\xi})^*= \alpha \ I \widetilde{\eta}_{\xi}$ for some
$\alpha \in \widetilde{F}.$
Since $(I \widetilde{\eta}_{\xi})^q=I  \widetilde{\eta}_{\mathfrak{e}},$
and $(I \widetilde{\eta}_{\mathfrak{e}})^*=I \widetilde{\eta}_{\mathfrak{e}}$
($I  \widetilde{\eta}_{\mathfrak{e}}$ is the unit of $C'$) then we obtain $\alpha^q=1.$ We also deduce $\alpha^2=1$ from $I \widetilde{\eta}_{\xi}=((I \widetilde{\eta}_{\xi})^*)^*=\alpha^2 \ I  \widetilde{\eta}_{\xi}.$ If $q$ is an odd prime number then $\alpha=1.$
If $q=2$ or $q=4$ then $\alpha \in \{ -1, 1\}.$

Then for any $\theta \in G$ we have $\theta=\xi^m$ for some integer $m.$ Hence, we can assume that $(\mathcal{X}_{\theta} \  \widetilde{\eta}_{\theta})^*=(I \  \widetilde{\eta}_{\xi^m})^* \cdot (\mathcal{X}_{\theta} \   \widetilde{\eta}_{\mathfrak{e}})^* =$
$((I \  \widetilde{\eta}_{\xi})^*)^m \cdot (\mathcal{X}_{\theta}^{\mathfrak{t}} \   \widetilde{\eta}_{\mathfrak{e}}) = \alpha^m \mathcal{X}_{\theta}^{\mathfrak{t}} \cdot \widetilde{\eta}_{\theta},$
where $\mathfrak{t}$ denotes the transpose or symplectic involution
on the matrix algebra. Therefore we obtain the alternative 3 or 4.

Consider the last case $G=\mathbb{Z}/4 \mathbb{Z},$ and $H=\{ \bar{0}, \bar{2} \}.$
Any element $a \in C=M_k(\widetilde{F}[H])$ can be uniquely represented in the form
\begin{equation}\label{aform}
a=\sum_{i=0}^{3} a_i =(a_0 + a_1) + (a'_0 +a'_1) (I \eta_{\bar{2}}),
\end{equation}
where $a_i \in C_{\bar{i}}$ are $G$-homogeneous components of $a,$
$a'_0=a_2 \eta_{\bar{2}} \in C_{\bar{0}},$ \  $a'_1=a_3 \eta_{\bar{2}} \in C_{\bar{1}}.$
Similarly to the previous case we obtain that $a^*=(a_0 + a_1)^* + (I \ \eta_{\bar{2}})^* (a'_0 +a'_1)^*=(a_0 +a_1)^* + \alpha (a'_0 +a'_1)^* \ (I \eta_{\bar{2}}),$ \ where $\alpha \in \{ -1, 1\},$ \ $I$ is the identity matrix. Hence it is enough to describe the restriction of our involution on the graded subspace $C_{\bar{0}} \oplus C_{\bar{1}}$ of $C.$

Let us take the vector space $A=C_{\bar{0}} \oplus C_{\bar{1}}.$ Define in $A$ the multiplication by the rule
$(a_0+a_1) \odot (b_0+b_1)=a_0 b_0 +a_0 b_1 + a_1 b_0 + a_1 b_1 \eta_{\bar{2}}$ in $C,$
where $a_i, b_i \in C_{\bar{i}}.$  It is clear that $A$ is a superalgebra ($(\mathbb{Z}/2 \mathbb{Z})$-graded algebra) with the
$(\mathbb{Z}/2 \mathbb{Z})$-grading $A=C_{\bar{0}} \oplus C_{\bar{1}}.$
The map $\bar{*}$ is naturally defined in $A$ by $(a_0+a_1)^{\bar{*}}=(a_0+a_1)^{*}=a_0^*+a_1^*,$ \ $a_i \in C_{\bar{i}}.$  Since the involution $*$ in $C$ is graded then $\bar{*}$ is a $(\mathbb{Z}/2 \mathbb{Z})$-graded linear operator of the second order which satisfies
$(a_i \odot b_j)^{\bar{*}}=\alpha^{i \cdot j} \  (b_j^{\bar{*}} \odot a_i^{\bar{*}}),$
where $a_i, b_i \in C_{\bar{i}}.$
A linear operator with all mentioned properties is called $\alpha$-involution.
It is clear that a $(1)$-involution is a graded involution on the superalgebra $A,$
and a $(-1)$-involution is a superinvolution on the superalgebra $A$.
We denote by $\Phi(C)$ the superalgebra $A$ with the multiplication $\odot$
and the $\alpha$-involution $\bar{*}$ obtained of the $(\mathbb{Z}/4 \mathbb{Z})$-graded algebra $C$
with a graded involution $*$ by the represented procedure.

The superalgebra $\Phi(C)$ is isomorphic as a superalgebra to the full matrix algebra $\mathcal{A}=M_k(\widetilde{F})$ with the elementary $(\mathbb{Z}/2 \mathbb{Z})$-grading defined by the $k$-tuple
$(\bar{\theta}_1, \dots, \bar{\theta}_k).$ Where $\bar{\theta}_i=\theta_i+H \in
(\mathbb{Z}/4 \mathbb{Z})/H$ if $(\theta_1,\dots,\theta_k)$ is the $k$-tuple defining
the $(\mathbb{Z}/4 \mathbb{Z})$-grading of $C.$ The graded isomorphism $\varphi:\Phi(C) \rightarrow M_k(\widetilde{F})$ is given by the rule $\varphi(E_{i j} \eta_{\xi_{i j}})=E_{i j}.$
Denote by $\tilde{*}$ the $\alpha$-involution on the superalgebra $\mathcal{A}$ induced
by the $\alpha$-involution $\bar{*}$ with respect to the isomorphism $\varphi.$
Hence we have $\varphi(a^{\bar{*}})=\varphi(a)^{\tilde{*}}$ for any $a \in \Phi(C),$ and
$\varphi$ is the isomorphism of superalgebras with $\alpha$-involution.

If $\alpha=1$ then $(\mathcal{A},\tilde{*})$ is isomorphic as a $*$-graded algebra to the matrix algebra $M_k(\widetilde{F})$ with an elementary grading and an elementary involution by \cite{BahtShestZ}, \cite{BahtZaic3}.
For $\alpha=-1$ superinvolutions on $(\mathbb{Z}/2 \mathbb{Z})$-graded matrix algebras are described by
M.L. Racine \cite{Racine} (Proposition 13, 14) (see also \cite{BahtTT}, Proposition 1).
In both of the cases $(\mathcal{A},\tilde{*})$ is isomorphic as a superalgebra with $\alpha$-involution
to the algebra $\widetilde{\mathcal{A}}=M_k(\widetilde{F})$ with an elementary grading
and an elementary $\alpha$-involution $\tilde{\star}.$ An $\alpha$-involution $\tilde{\star}$ on a matrix superalgebra $M_k(\widetilde{F})$ is called elementary if $(E_{i j})^{\tilde{\star}}=\pm E_{s t}$ for all $i,j=1,\dots,k$ and some $s, t.$

Consider the $k$-tuple $(\bar{\vartheta}_1,\dots,\bar{\vartheta}_k) \in \mathbb{Z}/2 \mathbb{Z}$
defining the elementary grading of $\widetilde{\mathcal{A}}.$ Suppose that
$\bar{\vartheta}_i=\vartheta_i + 2 \mathbb{Z},$ \ $\vartheta_i \in \{ 0, 1\},$
$i=1,\dots,k.$ Let us take the $k$-tuple $(\tilde{\vartheta}_1,\dots,\tilde{\vartheta}_k),$ where
$\tilde{\vartheta}_i=\vartheta_i + 4 \mathbb{Z} \in \mathbb{Z}/4 \mathbb{Z}.$
Denote by $\widetilde{C}$ the algebra $M_k(\widetilde{F}[H])$ with the canonical $(\mathbb{Z}/4 \mathbb{Z})$-grading defined by the $k$-tuple $(\tilde{\vartheta}_1,\dots,\tilde{\vartheta}_k)$ (see Lemma \ref{simple}).
Observe that $\Phi(\widetilde{C})$ is isomorphic as a superalgebra to $\widetilde{\mathcal{A}}$
by the arguments above. Suppose that $\tilde{\varphi}:\Phi(\widetilde{C}) \rightarrow \widetilde{\mathcal{A}}$ is the $(\mathbb{Z}/2 \mathbb{Z})$-graded isomorphism given by the rules
$\tilde{\varphi}(E_{i j} \eta_{\xi_{i j}})=E_{i j},$ \ $i, j=1,\dots,k.$

Define a map $\star$ on the algebra $\widetilde{C}$ by the rule:
\begin{equation} \label{star}
((c_0 + c_1) + (c'_0 +c'_1) (I \eta_{\bar{2}}))^\star=
\tilde{\varphi}^{-1}((\tilde{\varphi}(c_0 + c_1))^{\tilde{\star}}) +
\alpha \tilde{\varphi}^{-1}((\tilde{\varphi}(c'_0 + c'_1))^{\tilde{\star}}) (I \eta_{\bar{2}}).
\end{equation}
Where $c=(c_0 + c_1)  + (c'_0 +c'_1) (I \eta_{\bar{2}})$ is the form (\ref{aform}) for an element
$c \in \widetilde{C},$ \ $c_i,$ $c'_i \in \widetilde{C}_{\bar{i}}.$
It is clear that $\star$ is a $(\mathbb{Z}/4 \mathbb{Z})$-graded linear operator of $\widetilde{C}$
satisfying $(I \eta_{\bar{2}})^{\star}=\alpha (I \eta_{\bar{2}}).$
The restriction of $\star$ on $\Phi(\widetilde{C})$ is $\bar{\star}=\tilde{\varphi}^{-1} \tilde{\star} \tilde{\varphi}.$ Hence $\bar{\star}$ is an $\alpha$-involution, and $\tilde{\varphi}$ is an isomorphism
of superalgebras with $\alpha$-involution $(\Phi(\widetilde{C}),\bar{\star})$ and
$(\widetilde{\mathcal{A}},\tilde{\star}).$ Particularly, $(\Phi(\widetilde{C}),\bar{\star})$ is
isomorphic as a superalgebra with $\alpha$-involution to $(\Phi(C),\bar{*}).$
Moreover, for a basic element $b=E_{i j} \eta_{\xi}$ of the canonical basis of $\widetilde{C}$
we have
\begin{eqnarray} \label{star2}
&&b^{\star}=(E_{i j} \eta_{\xi})^{\star}=\alpha^{s} \cdot  E^{\tilde{\star}}_{i j} \  \eta_{\tilde{\xi}}, \\
&&s=\chi(\deg_{(\mathbb{Z}/4 \mathbb{Z})} b), \qquad \tilde{\xi}=\xi + \deg_{(\mathbb{Z}/4 \mathbb{Z})} E_{ij} - \deg_{(\mathbb{Z}/4 \mathbb{Z})} E^{\tilde{\star}}_{i j}. \nonumber
\end{eqnarray}
Using (\ref{star}) or (\ref{star2}) and Lemma \ref{Propchi} it can be directly checked that
$\star$ is a $(\mathbb{Z}/4 \mathbb{Z})$-graded involution of $\widetilde{C}.$
The isomorphism of superalgebras with $\alpha$-involution $\psi:(\Phi(\widetilde{C}),\bar{\star}) \rightarrow (\Phi(C),\bar{*})$ induces the isomorphism of $(\mathbb{Z}/4 \mathbb{Z})$-graded algebras
with involution $\Psi:(\widetilde{C},\star) \rightarrow (C,*)$ by the rule based on the representation (\ref{aform}) of elements $c \in \widetilde{C}$:
\[\Psi(c)=\Psi\bigl((c_0 + c_1)  + (c'_0 +c'_1) (I \eta_{\bar{2}})\bigr)=
\psi(c_0 + c_1)  + \psi(c'_0 +c'_1) (I \eta_{\bar{2}}).\]
Since $\star$ is an elementary involution (by (\ref{star2})) then we obtain the last alternative of the theorem. It completes the proof.
\hfill $\Box$

The direct consequence of Theorem \ref{simGr-pClas} is the next corollary.
\begin{corollary}  \label{assum-cicl}
Let $\widetilde{F}$ be an algebraically closed field of
characteristic zero.  Assumption \ref{Class} is true
over $\widetilde{F}$ for a cyclic group $G$ of a prime order or
of the order $4.$
\end{corollary}

Hence Theorems \ref{simGr-pClas}, \ref{*PI-gen2}
immediately imply Theorem \ref{*PI}.
\begin{theorem}  \label{*PI}
Let $q$ be a prime number or $q=4,$ $F$ a field of
characteristic zero. Then for any $(\mathbb{Z}/q \mathbb{Z})$-graded
finitely generated associative PI-algebra $A$ with graded involution over $F$
there exists a finite dimensional over $F$ $(\mathbb{Z}/q \mathbb{Z})$-graded associative algebra $C$
with graded involution such that the ideals of graded identities with involution of $A$ and $C$ coincide.
\end{theorem}

It is an interesting problem to describe all groups $G$ such that
Assumption \ref{Class} is true for $G$-graded $*$-algebras over an algebraically closed field.
We suppose that Assumption \ref{Class} should be true for any
finite abelian group.

\begin{conjecture} \label{ClFinAb}
Let $G$ be a finite abelian group, and
$\widetilde{F}$ an algebraically closed field of characteristic zero.
Given a $G$-graded finite dimensional algebra $A$ with graded involution
$A$ is $*$-graded simple if and only if $A$ is isomorphic as a
graded $*$-algebra either to $G$-graded simple algebra
$\widetilde{C}^{(1)}=M_k(\widetilde{F}^{\zeta}[H])$
with an elementary involution, or to the direct product
$\widetilde{C}^{(2)}=\mathcal{B} \times \mathcal{B}^{op}$ of
a graded simple algebra $\mathcal{B}=M_k(\widetilde{F}^{\zeta}[H])$
and its opposite algebra $\mathcal{B}^{op}$ with the exchange involution $\bar{*}.$
Where $H$ is a subgroup of $G,$ and
$\zeta:H \times H \rightarrow \mathbb{Q}[\sqrt[\mathfrak{m}]{1}]^{*}$
is a 2-cocycle on $H$ with values in the algebraic extension of rational
numbers $\mathbb{Q}$ by a primitive root $\sqrt[\mathfrak{m}]{1}$ of $1$
of degree $\mathfrak{m}=|G|.$
\end{conjecture}

If it is true then any $G$-graded finitely generated PI-algebra with graded involution
over a field of characteristic zero
should be PI-representable with respect to graded $*$-identities.

Both of the questions (the classification of finite dimensional $*$-graded simple algebras,
and PI-representability of finitely generated algebras) are also interesting in case
of a finite (not necessary abelian) group.

\section{Acknowledgments.}
The work was partially supported  by CNPq, CAPES, and CNPq-FAPDF PRONEX
grant 2009/00091-0 (193.000.580/2009).
The author is specially thankful to FAPESP for the financial support and possibility
to finish this work, and personally deeply appreciate to Ivan Shestakov, Antonio
Giambruno and Victor Petrogradsky for constant inspiration and
useful discussions.

\bibliographystyle{amsplain}

\end{document}